\documentclass[aap,citesort,MSNbibl,dvips]{arximspdf}
\usepackage{dcolumn}
\usepackage{graphicx}
\doi{10.1214/10-AAP705}
\volume{21}
\issue{3}
\pubyear{2011}
\firstpage{876}
\lastpage{907}

\makeatletter

\newcolumntype{d}[1]{D{.}{.}{#1}}

\newcommand{\bigtimes}{\mathop{\mbox{\fontsize{17}{17}\selectfont{$\!\times$}}}}

\newcommand{\mb}[1]{\mathbf{#1}}
\newcommand{\mbb}[1]{\mathbb{#1}}

\newtheorem{theorem}{Theorem}
\newtheorem{lem}[theorem]{Lemma}
\newtheorem{prop}[theorem]{Proposition}
\newproclaim{defn}{Definition}
\newproclaim{assump}{Assumption}
\newproclaim{example}{Example}
\newproclaim{Remark}{Remark}

\makeatother

\begin{document}
\begin{frontmatter}

\title{General notions of indexability for queueing control and asset management}
\runtitle{General notions of indexability}

\begin{aug}
\author[A]{\fnms{Kevin D.} \snm{Glazebrook}\corref{}\thanksref{t1}\ead[label=e1]{k.glazebrook@lancaster.ac.uk}\ead[label=u1,url]{http://www.lums.lancs.ac.uk/profiles/kevin-glazebrook/}},
\author[B]{\fnms{David J.} \snm{Hodge}\thanksref{t1}\ead[label=e2]{david.hodge@nottingham.ac.uk}} and
\author[C]{\fnms{Chris}~\snm{Kirkbride}\thanksref{t2}\ead[label=e3]{c.kirkbride@lancaster.ac.uk}}
\runauthor{K. D. Glazebrook, D. J. Hodge and C. Kirkbride}
\affiliation{Lancaster University}
\address[A]{K. D. Glazebrook \\
Department of Mathematics and Statistics \\
and\\
Department of Management Science\\
Lancaster University\\
LA1 4YF\\
United Kingdom\\
\printead{e1}\\
\printead{u1}}
\address[B]{D. J. Hodge \\
School of Mathematical Sciences\\
University of Nottingham\\
NG7 2RD\\
United Kingdom\\
\printead{e2}}
\address[C]{C. Kirkbride\\
Department of Management Science\\
Lancaster University Management School\\
LA1 4YX\\
United Kingdom\\
\printead{e3}}
\end{aug}

\thankstext{t1}{Supported by EPSRC Grant EP/E049265/01.}
\thankstext{t2}{Supported by an RCUK Fellowship.}

\received{\smonth{8} \syear{2009}}
\revised{\smonth{3} \syear{2010}}

%
\begin{abstract}
We develop appropriately generalized notions of \textit{indexability} for
problems of dynamic resource allocation where the resource concerned
may be assigned more flexibility than is allowed, for example, in
classical multi-armed bandits. Most especially we have in mind the
allocation of a divisible resource (manpower, money, equipment) to a
collection of objects (projects) requiring it in cases where its
over-concentration would usually be far from optimal. The resulting
\textit{project indices} are functions of both a resource level and a
state. They have a simple interpretation as a \textit{fair charge} for
increasing the resource available to the project from the specified
resource level when in the specified state. We illustrate ideas by
reference to two model classes which are of independent interest. In
the first, a pool of servers is assigned dynamically to a collection of
service teams, each of which mans a service station. We demonstrate
indexability under a natural assumption that the service rate delivered
is increasing and concave in the team size. The second model class is a
generalization of the \textit{spinning plates model} for the optimal
deployment of a divisible investment resource to a collection of reward
generating assets. Asset indexability is established under
appropriately drawn laws of diminishing returns for resource
deployment. For both model classes numerical studies provide evidence
that the proposed \textit{greedy index heuristic} performs strongly.
\end{abstract}

%
\begin{keyword}[class=AMS]
\kwd[Primary ]{68M20}
\kwd[; secondary ]{90B22}
\kwd{90B36}.
\end{keyword}
\begin{keyword}
\kwd{Asset management}
\kwd{dynamic programming}
\kwd{dynamic resource allocation}
\kwd{full indexability}
\kwd{index policy}
\kwd{Lagrangian relaxation}
\kwd{monotone policy}
\kwd{queueing control}.
\end{keyword}

\end{frontmatter}

\section{Introduction}
A notable, now classical, contribution to the theory of dynamic
resource allocation was the elucidation by Gittins \cite{git79,git89}
of \textit{index-based solutions} to a large family of \textit{multi-armed
bandit problems} (MABs). This is a class of models concerned with the
sequential allocation of \textit{effort}, to be thought of as a \textit
{single indivisible resource}, to a collection of \textit{stochastic
reward generating projects} (or \textit{bandits} as they are sometimes
called). Gittins demonstrated that optimal project choices are those of
\textit{highest index}. There is no doubt that the idea that strongly
performing policies are determined by simple, interpretable
calibrations (i.e., \textit{indices}) of decision options is an
attractive and powerful one and offers crucial computational benefits.
There is now substantial literature describing extensions to and
reformulations of Gittins' result. Some key contributions are cited in
the recent survey of Mahajan and Teneketzis \cite{mah07}.

Whittle \cite{whi88} introduced a class of \textit{restless bandit
problems} (RBPs) as a means of addressing a critical limitation of
Gittins' MABs, namely, that projects should remain frozen while not in
receipt of effort. In RBPs, projects may change state while active or
passive though according to different dynamics. However, this
generalization is bought at great cost. In contrast to MABs, RBPs are
almost certainly intractable having been shown to be \textit{PSPACE-hard}
by Papadimitriou and Tsitsiklis~\cite{pap99}. Whittle \cite{whi88}
proposed an index heuristic for those RBPs which pass an \textit
{indexability test}. This heuristic reduces to Gittins' index policy in
the MAB case. Whittle's index emerges from a Lagrangian relaxation of
the original problem and has an interpretation as a \textit{fair charge}
for the allocation of effort to a particular project in a particular
state. Weber and Weiss \cite{web90} established a form of asymptotic
optimality for Whittle's heuristic under given conditions. More
recently, several studies have demonstrated the power of Whittle's
approach in a range of application areas. These include the dynamic
routing of customers for service \cite{arg08,gla07}, machine
maintenance \cite{gla05}, asset management \cite{gla06} and inventory
routing \cite{arc08}.

The above classical models and associated theory are undeniably
powerful when applicable. However, the scope of their applicability is
heavily constrained by the very simple view the models take of the
resource to be allocated. As indicated above, in Gittins' MAB model a
single indivisible resource is allocated wholly and exclusively to a
single project at each decision epoch. In Whittle's RBP formulation,
parallel server versions of this are allowed. Many applications,
however, call for the allocation of a \textit{divisible resource} (e.g.,
money, manpower or equipment) in situations where its over
concentration would usually be far from optimal. This is the case, for
example, in the problem concerning the planning of new product
pharmaceutical research which was discussed by Gittins \cite{git89} and
which provided practical motivation for his pioneering contribution.
This paper records the first outcomes of a major research program whose
goal is to develop a usable and effective index theory for such problems.

In Section~\ref{sec2} we present a general model for dynamic resource
allocation. Both Gittins' MABs and Whittle's RBPs may be recovered as
special cases as may the recent model of\vadjust{\goodbreak} \cite{gla08} which extends
Gittins' MABs such that bandit activation consumes amounts of the
available resource which may vary by bandit and state. Our general
model allows for resource to be applied at a range of levels to each
constituent project, subject to some overall constraint on the total
rate at which resource is available. A notion of (\textit{full}) \textit
{indexability} which generalizes that of Whittle for RBPs is developed.
Any project which is \textit{fully indexable} has an index which is a
function \textit{both} of a given resource level ($a$) and of a given
state ($x$). The index $W(a,x)$ may be understood as a \textit{fair
charge} for raising the project's resource level above $a$ when in
state $x$. We discuss how to use such indices to develop heuristics for
dynamic resource allocation when all projects are fully indexable.

In Sections~\ref{sec3} and~\ref{sec4} we use the ideas and methods of Section~\ref{sec2} to
construct index heuristics for the dynamic allocation of a divisible
resource in the context of two model classes which are of considerable
interest in their own right. In Section~\ref{sec3} we deploy the framework of
Section~\ref{sec2} to develop heuristics for the dynamic allocation of a \textit
{pool of $S$ servers} to $K$ service stations (or customer classes) at
which queues may form. This model is able to capture situations where,
for example, each of $K$ customer classes is served by a dedicated team
of specialists. Additionally, $S$ higher level generalist servers are
available for deployment across the customer classes to supplement the
specialist teams as demand dictates. Deployment of $a_k$ generalists to
customer class $k$ enhances the local specialist team which then
delivers service collectively at rate $\mu_k(a_k)$. An assumption that
the \textit{service rate} functions $\mu_k$ are increasing and concave
reflects a law of diminishing returns as service teams grow. The
problem of determining how the pool of generalists should be deployed
across the customer classes in response to queue length information is
formulated as a dynamic resource allocation problem of the kind
discussed in Section~\ref{sec2}. The analysis which establishes full
indexability in Section~\ref{sec3} markedly adds to the queueing control
literature in establishing monotonicity with respect to service costs
of optimal policies for a derived problem involving a single queue. An
algorithm is given for the computation of indices. A numerical study
provides evidence that a \textit{greedy index heuristic} for allocating
the common service pool is close to optimal throughout a numerical
study featuring nearly 10,000 two station problems.

The model class studied in Section~\ref{sec4} generalizes the so-called
\textit {spinning plates model} discussed by Glazebrook, Kirkbride and
Ruiz-Hernandez \cite{gla06}. It is a flexible finite state model class
in which a divisible investment resource is available to drive
improvements to the (reward) performance of $K$ reward generating
assets, which in the absence of any such resource deployment will tend
to deteriorate. Positive investment \textit {both} arrests an asset's
tendency to deteriorate and enhances asset performance by enabling
movement of the asset state toward those in which its reward generating
performance will be stronger. Full indexability for assets is
established under laws of diminishing returns as asset investment
levels grow. This considerably extends the work of Glazebrook,
Kirkbride and Ruiz-Hernandez~\cite{gla06}. A numerical study which
features 14,000 two asset problems testifies to the strong performance
of the greedy index heuristic in comparison to optimum and to
competitor policies. Conclusions and proposals for further work are
discussed in Section~\ref{sec5}.

\section{A model for dynamic resource allocation}\label{sec2}

We propose a semi-Markov decision process (SMDP) formulation $\{
( \Omega_{k},L_{k},c_{k},r_{k},q_{k}), 1\leq k\leq$ $K\}$
of the problem of dynamically allocating a resource to a collection of $K$
stochastic projects. This formulation includes Gittins' MABs
and Whittle's RBPs as special cases. In our SMDP project $k$
is characterized by its (finite or countable) \textit{state space}
$\Omega
_{k}$, its \textit{highest activation level} $L_{k}\in\mathbb{Z}^{+}$,
\textit{cost rate function} $c_{k}\dvtx\{ 0,1,\ldots,L_{k}\}
\times\Omega_{k}\rightarrow\mathbb{R}^{+}$, \textit{resource
consumption function} $r_{k}\dvtx\{ 0,1,\ldots,\break L_{k}\} \times
\Omega_{k}\rightarrow\mathbb{R}^{+}$ and \textit{Markov transition law}
$q_{k}$. The model is in continuous time. We use $x_{k},x_{k}^{\prime}
\in\Omega_{k}$ for generic states of project $k$ and $\mathbf
{x},\mathbf{x}^{\prime} \in\bigtimes_{k=1}^{K} \Omega_{k}$ for
generic states of the process. In the
SMDP an action $\mathbf{a}=( a_{1},a_{2},\ldots,a_{K})$ must
be taken at time~$0$ and after each (state) transition of the process. This
specifies the \textit{resource level} $a_{k}\in\{ 0,1,\ldots
,L_{k}\}$ to be applied to project $k, 1\leq k\leq K$. The
choice $a_{k}=0$ indicates that resource at a minimal level (usually
none) is to be applied to $k$ ($k$ is \textit{passive}), while the choice
$a_{k}=L_{k}$ indicates a maximal resource
allocation. Resource level $a_{k}$ applied to project $k$ when in state
$x_{k}$ leads to a consumption of resource at rate
$r_{k}(a_{k},x_{k})$, with $r_{k}(\cdot,x_{k})$ increasing $\forall
k, x_{k}$. In the major
examples discussed in the upcoming sections we will have
$r_{k}(a_{k},x_{k})=a_{k}$ $\forall k,  x_{k}$ and the resource
level is identified with the resource consumed. When resource level
$a_k$ is applied to project $k$ when in state $x_k$, it incurs costs at
rate $c_k(a_k,x_k)$. Both cost and resource consumption rates are
additive over projects. It will be convenient to write $c(\mathbf
{a},\mathbf{x}) = \sum_k c_k(a_k,x_k)$ and $r(\mathbf{a},\mathbf
{x}) =
\sum_k r_k(a_k,x_k)$. The set of \textit{admissible actions} in process
state $\mathbf{x}$ is given by $A(\mathbf{x})=\{ \mathbf{a};
r(\mathbf
{a},\mathbf{x}) \leq R \}$ where $R$ is the rate at which resource is
available to the system, assumed constant over
time. We suppose that $A(\mathbf{x})\neq\phi,  \mathbf{x} \in
\bigtimes_{k=1}^{K} \Omega_{k}$. An \textit{admissible policy} is a rule for
taking admissible actions.

Should action $\mathbf{a}$ be taken when the system is in state
$\mathbf{x}$,
the system will remain in state $\mathbf{x}$ for an amount of time
which is
exponentially distributed with rate
\[
\sum_{\mathbf{x}^{\prime} \in\bigtimes_{k} \Omega_{k}} \mathbf{q}
(\mathbf{x}^{\prime}\mid\mathbf{x},\mathbf{a}) = \sum_{k=1}^{K}
\sum_{x_{k}^{\prime} \in\Omega_{k}} q_{k}(x_{k}^{\prime} \mid
x_{k},a_{k}) \leq Q < \infty\qquad \forall  \mathbf{x},  \mathbf{a}.
\]
The transition following will be from state $x_{k}$ to state
$x_{k}^{\prime}$ within project $k$ with probability
\[
q_{k}(x_{k}^{\prime} \mid x_{k},a_{k}) \biggl\{ \sum_{\mathbf
{x}^{\prime
} \in\,\bigtimes_{k} \mathbf{\Omega}_{k}}\mathbf{q}(\mathbf
{x}^{\prime}
\mid\mathbf{x},\mathbf{a}) \biggr\}^{-1}.
\]
Hence the projects evolve independently, given the choice of action,
with $q_{k}$ yielding transition rates for project $k$. The goal of
analysis is
the determination of a policy for resource allocation (a rule for taking
admissible actions at all decision epochs) which minimizes the average cost
per unit time incurred over an infinite horizon.

To develop ideas and notation we use $\bar{\mathbf{U}}$ for the set of
\textit{deterministic}, \textit{stationary}, \textit{Markov} (\textit{DSM}) \textit{and admissible policies}
determined by functions $\mb{u}$ with domain $\bigtimes_{k=1}^K
\Omega
_k$ which satisfy $\mb{u}(\mb{x}) \in A(\mb{x})$  $\forall \mb{x}$.
Fix $\mb{u} \in\bar{\mb{U}}$. We shall also use $\{\mb{X}(t), t
\geq
0\}$ for the system state evolving over time and $[\mb{u}\{\mb{X}(t)\},
t \geq0]$ for the corresponding stochastic process of admissible
actions taken by $\mb{u}$. We write
%
\begin{equation} \label{n1}
C(\mb{u},\mb{x}) = \mathop{\lim\inf}_{t \rightarrow\infty}
\frac{1}{t}
\biggl( \int_0^t \mbb{E}_{\mb{u}}^{\mb{x}} c ( \mb{u}
\{ \mb
{X}(s) \}, \mb{X}(s) ) \,ds \biggr)
\end{equation}
for the average cost per unit time incurred under policy $\mb{u}$ over
an infinite horizon from initial state $\mb{x}$. In~(\ref{n1}) $\mbb
{E}_{\mb{u}}^{\mb{x}}$ denotes an expectation taken over realizations
of the system evolving under $\mb{u}$ from initial state $\mb{x}$. We
shall assume the existence of a policy $\mb{u} \in\bar{\mb{U}}$ such
that $C(\mb{u},\mb{x}) < \infty$ $\forall \mb{x}$ and write
$C^{\mathrm{opt}}(\mb{x})$ for the minimized cost rate, namely,
%
\begin{equation} \label{n2}
C^{\mathrm{opt}}(\mb{x}) = \inf_{\mb{u} \in\bar{\mb{U}}} C(\mb{u},\mb{x}).
\end{equation}
We shall use the term \textit{optimal} to denote a policy (assumed to
exist) which achieves the infimum in~(\ref{n2}) uniformly over initial
states. This applies both to the problem in~(\ref{n2}) and also to the
derived optimization problems we shall discuss later in the account. In
the model classes featured in Sections~\ref{sec3} and~\ref{sec4} it will be the case that
the average costs in~(\ref{n1}) and~(\ref{n2}) are independent of
$\mb
{x}$. Henceforth, for simplicity, we shall suppress dependence on the
initial state $\mb{x}$ in the notation.

We shall use
%
\begin{equation} \label{n3}
R(\mb{u}) = \mathop{\lim\inf}_{t \rightarrow\infty}  \frac
{1}{t}  \biggl(
\int_0^t \mbb{E}_{\mathbf{u}} r(\mb{u}\{\mb{X}(s)\},\mb{X}(s)) \,ds \biggr)
\end{equation}
for the average rate at which resource is consumed under policy $\mb
{u}$. We also write
%
\begin{equation} \label{n4}
C(\mb{u}) = \sum_{k=1}^K C_k(\mb{u}), \qquad R(\mb{u}) = \sum_{k=1}^K
R_k(\mb{u})
\end{equation}
to give a disaggregation of the cost and resource consumption rates
into the contributions from individual projects.

In principle, the tools of dynamic programming (DP) are available to
determine optimal policies. See, for example, \cite{put94}. However,
direct application of DP is
computationally infeasible other than for small problems (crucially,
small $K$). Hence, our primary interest lies in the development of
\textit{heuristic policies} which are close to cost minimizing. To this end we
relax the
optimization problem in~(\ref{n2}) by extending the class of policies
from the DSM admissible class $\bar{\mb{U}}$ to those DSM policies
$\mathbf{u}\dvtx \bigtimes_{k=1}^{K} \Omega_{k} \rightarrow
\bigtimes_{k=1}^{K} \{0,1,\ldots,L_{k}\}$ which consume resource
at an \textit{average rate} which is no greater than $R$. Hence, we write
%
\begin{equation} \label{n5}
\acute{C}^{\mathrm{opt}}=\inf_{\mathbf{u}}\sum_{k=1}^{K}C_{k}(\mathbf{u}),
\end{equation}
where in~(\ref{n5}), the infimum is taken over the collection of
DSM policies satisfying
%
\begin{equation} \label{n6}
\sum_{k=1}^{K}R_{k}(\mathbf{u})\leq R.
\end{equation}
We now relax the problem again by further extending the class of policies
and by incorporating the constraint~(\ref{n6}) into the objective
(\ref{n5}) in a Lagrangian fashion. We write
%
\begin{equation} \label{n7}
C(W)=\inf_{\mathbf{u}}\sum_{k=1}^{K} \{C_{k}(\mathbf
{u})+WR_{k}(\mathbf
{u})\}-WR.
\end{equation}
In~(\ref{n7}) the infimum is taken over the class of DSM policies
$\mathbf{u}\dvtx\bigtimes_{k=1}^{K} \Omega_{k} \rightarrow
\bigtimes_{k=1}^{K} \{0,1,\ldots,L_{k}\}$ which allow, for each project
$k$, a free choice
of action from the set $\{ 0,1,\ldots, L_{k}\} $ at each
decision epoch. It is clear that
\[
C(W) \leq\acute{C}^{\mathrm{opt}} \leq C^{\mathrm{opt}},\qquad  W \in\mathbb{R}^{+}.
\]
However, the Lagrangian relaxation of our optimization problem
expressed by
(\ref{n7}) admits, on account both of the policy class involved and the
nature of the objective, an additive project-based decomposition. Expressed
differently, an optimal policy for~(\ref{n7}) operates optimal policies
for the individual projects in parallel. In an obvious notation we write
%
\begin{equation} \label{e5}
C(W) = \sum_{k=1}^{K} C_{k}(W)-WR,
\end{equation}
where
%
\begin{equation} \label{e6}
C_{k}(W) = \inf_{u_{k}} \{ C_{k}(u_{k})+WR_{k}(u_{k})\},\qquad
1 \leq k \leq K.
\end{equation}
The optimization problem in~(\ref{e6}) concerns \textit{project $k$
alone}. We denote it $P(k,W)$. In its objective the Lagrange multiplier
$W$ plays the role of a charge per unit of time and per unit of
resource consumed. An optimal policy $u_{k}(W)$ for $P(k,W)$ minimizes
an aggregate rate of project costs incurred and charges levied for
resource consumed. Further, the policy $\mathbf{u}(W)$ which applies
$u_{k}(W)$ to each project~$k$, achieves $C(W)$ in (\ref {n7}) and
hence provides a solution to the above Lagrangian relaxation. Note that
in what follows we shall use the notation $\mathbf {u}(W,\mathbf{x}),
u_{k}(W,x_{k})$ to denote the action (resource consumption levels)
chosen by DSM policies $\mathbf{u}(W), u_{k}(W)$ in states $\mathbf{x},
x_{k}$, respectively.

In order to develop natural \textit{project calibrations} (or \textit
{indices}) which can facilitate the construction of effective
heuristics for
our original problem~(\ref{n2}), we seek optimal policies for the
problems $\{ P(k,W),W \in\mathbb{R}^{+}, 1 \leq k \leq K
\}
$ which are structured as in Definition~\ref{defn1} below.
We first require additional notation. Write
%
\begin{equation} \label{e7}\quad
\Pi_{k} \{u_{k}(W),a\} = \{ x \in\Omega_{k}; u_{k}(W,x) \leq a
\},\qquad a \in\{ 0,1,\ldots, L_{k}-1\},
\end{equation}
for the set of project $k$ states for which policy $u_{k}(W)$ chooses to
consume resource at level $a$ or below.
\begin{defn}[(Full indexability)] \label{defn1}
Project $k$ is fully indexable if there exists a family of DSM policies
$\{ u_{k}(W), W \in\mathbb {R}^{+}\} $ such that $u_k(W)$ is optimal
for $P(k,W)$  $\forall W$ and $\Pi_{k} \{ u_{k}(W),a\}$ is
nondecreasing in $W$ for each $a \in\{ 0,1,\ldots,L_{k}-1 \}$.
\end{defn}

To summarize the requirements of Definition~\ref{defn1}, a project $k$
will be fully indexable if the problem $P(k,W)$ has an optimal policy
which, for any given state, consumes an amount of resource which is
\textit{decreasing} in the resource charge $W$. Full indexability
enables a \textit{calibration} of the individual projects as described
in Definition~\ref{defn2}.
\begin{defn}[(Project indices)] \label{defn2}
If project $k$ is fully indexable as in
Definition~\ref{defn1}, a corresponding index function $W_{k}\dvtx \{
0,1,\ldots,L_{k}-1\} \times\Omega_{k} \rightarrow\mathbb
{R}^{+}$ is given by
%
\begin{equation} \label{e8}
W_{k}(a,x) = \inf[W; x \in\Pi_{k} \{u_{k}(W),a\}].
\end{equation}
\end{defn}
\begin{Remark*}
The index $W_{k}(a,x) $ can be thought of as a \textit{fair charge} at
project $k$ for
raising the resource level from $a$ to $a+1$ in state $x$. Were a resource
charge less than $W_{k}(a,x)$ to be levied, the consumption of the additional
resource would be preferable, while if the resource charge were to be in
excess of the index, that would not be the case. We shall adopt the
convention that the index function is extended to $W_{k}\dvtx \{
-1,0,1,\ldots,L_{k}\} \times\Omega_{k} \rightarrow\mathbb
{R}^{+} \cup\{ \infty\}$ where $W_{k}(-1,x) = \infty,
W_{k}(L_{k},x)=0$  $\forall x \in\Omega_{k}$.
\end{Remark*}

The following is a simple consequence of the above definitions. Its
proof is omitted.\vadjust{\goodbreak}
\begin{lem} \label{lem1}
If project $k$ is fully indexable, the index $W_{k}(a,x)$ is decreasing
in~$a$, for fixed $x$.
\end{lem}

Hence, under full indexability, the fair charge for raising the resource
level for project $k$ in any state $x$ from $a$ to $a+1$ is decreasing in
the resource level $a$.

We now return to consideration of the Lagrangian relaxation in (\ref
{n7}) and~(\ref{e5})
and suppose that all $K$ projects are fully indexable with families of
optimal policies
\[
\{ u_{k}(W), W \in\mathbb{R}^{+}, 1 \leq k \leq K\}
\]
structured as in Definition~\ref{defn1}. Under full indexability, all
of these
policies have a structure describable in terms of the index functions
$W_{k}, 1 \leq k \leq K$. Theorem~\ref{thm2} now follows.
\begin{theorem} \label{thm2}
Suppose that all $K$ projects are fully indexable with extended index
functions $W_{k}\dvtx \{ -1,0,1,\ldots,L_{k} \} \times\Omega
_{k} \rightarrow\mathbb{R}^{+} \cup\{ \infty\}$. The
policy $\mathbf{u}(W)$ such that
\begin{eqnarray}
\mathbf{u}(W,\mathbf{x}) = \mathbf{a} \quad\Longleftrightarrow\quad
W_{k}(a_{k}-1,x_{k}) > W \geq W_{k}(a_{k},x_{k}),\nonumber\\
&&\eqntext{1\leq k\leq K, \mathbf{x} \in\bigtimes_{k=1}^{K}
\Omega_{k},}
\end{eqnarray}
achieves $C(W)$ $\forall W \in\mbb{R}^+$.
\end{theorem}
\begin{Remark*}
According to Theorem~\ref{thm2}, policy $\mathbf{u}(W)$ constructs
actions (allocations of resource) in each system state by
accumulating resource at each project until the fair charge for adding
further resource drops below the prevailing charge $W$. This is strongly
suggestive of how effective, interpretable heuristics for our original
dynamic resource allocation problem based on the above indices (fair
charges) may be constructed when all projects are fully indexable. A natural
\textit{greedy index heuristic} constructs actions in every system state
by increasing resource consumption levels in \textit{decreasing order} of
the above station indices until the point is reached when the resource
constraint is violated by additional allocation of resource.
\end{Remark*}

Formally the greedy index heuristic is structured as follows:

\subsection*{Greedy index heuristic} In state $\mathbf{x}$ the
greedy index heuristic
constructs an action (allocation of resource) as follows:

\textit{Step} 1. The initial allocation is $\mathbf
{0}=\{ 0,0,\ldots,0\}$. The current allocation is
$\mathbf
{a}=\{ a_{1},a_{2},\ldots,a_{K}\}$ with $\sum
_{k}r_{k}(a_{k},x_{k})<R$.

\textit{Step} 2. Choose any $k$ satisfying
\[
W_{k}(a_{k},x_{k})=\max_{1\leq j\leq K}W_{j}(a_{j},x_{j}).
\]\vfill\eject

\textit{Step} 3. If $\mathbf{e}_k$ denotes a $K$-vector
whose $k$th component is 1 with zeroes elsewhere, the new deployment is
$\mathbf{a}+\mathbf{e}_{k}$ if
%
\begin{equation} \label{e9}
\sum_{l\neq k} r_{l}(a_{l},x_{l})+r_{k}(a_{k}+1,x_{k}) \leq R.
\end{equation}
If there is strict inequality in~(\ref{e9}), return to Step 1 and
repeat. Otherwise, stop and declare $\mathbf{a}+\mathbf{e}_{k}$ to be the
chosen action in $\mathbf{x}$. If
\[
\sum_{l\neq k}r_{l}(a_{l},x_{l})+r_{k}(a_{k}+1,x_{k})>R,
\]
stop and declare $\mathbf{a}$ to be the chosen action in $\mathbf{x}$.
\begin{Remark*}
We shall use Figure~\ref{fig1} to illustrate the construction of
actions by both the policy $\mathbf{u}(W)$ (as in Theorem~\ref{thm2})
and the greedy index heuristic in a simple problem with $K=2$ in which
both projects are fully indexable. Section~\ref{sec3} discusses a class of
models in which $r_k(a_k,x_k) = a_k$  $\forall k, x_k$ and where all
projects have state space $\mathbb{N}$ and a common maximum resource
level, $L$ say, which is equal to $R$, the total rate at which resource
is available. Suppose now that $L=R=5$ in such a model and that the
system state is $\mathbf{x} = (x_1,x_2) = (5,2)$. Figure~\ref{fig1}
indicates values of the appropriate project indices $W_1(a,5)$ and
$W_2(a,2)$ for the range $0 \leq a \leq4$ together with the value of
the Lagrange multiplier $W$.

\begin{figure}

\includegraphics{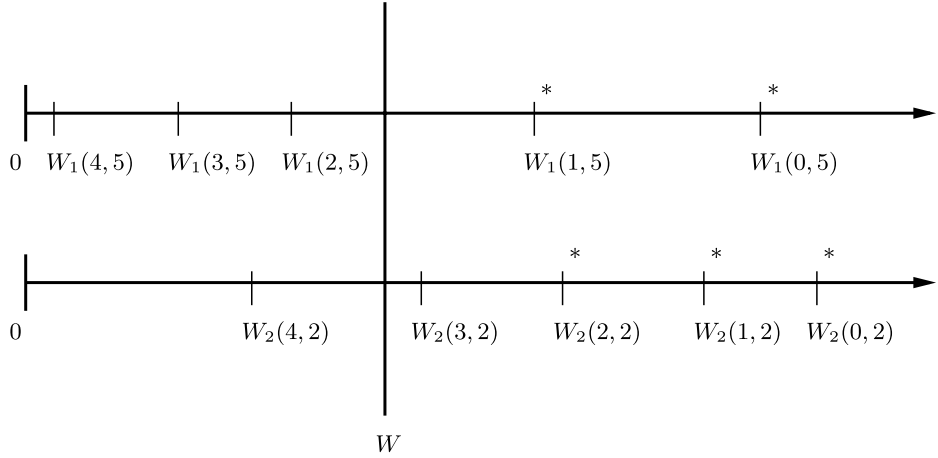}

\caption{Index values for state $\mathbf{x}=(5,2)$.} \label{fig1}
\end{figure}

The policy $\mathbf{u}(W)$ will make allocations of resource supported
by those index values which are above $W$. Hence from Figure \ref
{fig1}, the choice of action in state $\mathbf{x}=(5,2)$ will be
$\mathbf{a}=(2,4)$. This is an inadmissible action for the original
problem since the total resource rate allocated (6) exceeds that
available~(5).\vadjust{\goodbreak} The greedy heuristic makes allocations of resource
supported by the five largest index values (indicated by $*$ in Figure
\ref{fig1}). Plainly, the action taken by the index heuristic is
$\mathbf{a}=(2,3)$. As the system state evolves under the operation of
either policy, the index values change as do the implied actions.
\end{Remark*}

The major challenge to implementation of the above program for heuristic
construction is the identification of optimal policies for the problems
\[
\{ P(k,W),1\leq k\leq K,W\in\mathbb{R}^{+}\},
\]
which meet the requirements of Definition~\ref{defn1}. In Sections~\ref{sec3}
and~\ref{sec4} we are
able to achieve this in the context of two model classes for which
we are able to establish an appropriate form of full indexability. For
the Section~\ref{sec3} problem, we also give an algorithm for index computation.
For both model classes we proceed to assess the
performance of the greedy index heuristic in extensive numerical studies.

\begin{Remark*}
We recover Whittle's RBPs \cite{whi88} by making the choices
$r_{k}(a_{k}$, $x_{k})=a_{k}, L_{k}=1, 1\leq k\leq K$ and $R<K$ in
the above.
Hence there are just two modes of activation (active, passive) of each
project, with $R$ projects to be made active at each epoch. For this special
case the above greedy index heuristic is precisely the index heuristic proposed
by Whittle. If we make the further choice $R=1$ and impose the requirement
that projects can only change state under the active action, we then recover
Gittins' MAB \cite{git79} and its associated (optimal) index policy.
\end{Remark*}

\section{The optimal allocation of a pool of servers}\label{sec3}

We illustrate the above ideas by considering a set-up in which service
is provided at $K$ service stations.
These stations could represent distinct geographical locations or facilities
dedicated to the service of a particular class of customer. Customers arrive
at the stations in $K$ independent Poisson streams, with $\lambda_{k}$ the
rate for station $k$. A pool of $S$ servers is available to support service
at the $K$ stations. Should $a$ servers from the pool be allocated to
station $k$ at any point, the resulting exponential service rate is
$\mu
_{k}(a)$. Note that there may be a local team of servers
permanently stationed at $k$ (i.e., in addition to any allocated from the
pool) in which case we will have $\mu_{k}(0)>0$. Please note also that we
shall suppose that all servers (whether permanently based at a location or
allocated there from the common pool) offer service as a \textit{team},
namely, that they act in concert as a single server. The goal of
analysis is
the determination of a policy for deploying the common service pool in
response to queue length information to minimize some linear measure of
holding cost rate for the system incurred over an infinite horizon.

More formally, the \textit{system state} at time $t$ is $\mathbf
{n}(t)=\{ n_{1}(t),n_{2}(t),\ldots,n_{K}(t)\}$ where
$n_{k}(t)$ is
the number of customers at service station $k$ (including any in
service) at\vadjust{\goodbreak}
$t$. We shall on occasion refer to $n_{k}(t)$ as the \textit{head count}
at station $k$ at time $t$. This system state is observed continuously. The
\textit{decision epochs} for the system are time zero and the times at
which the system
state changes. At each decision epoch, some action $\mathbf{a}=(
a_{1},a_{2},\ldots,a_{K})$ is taken, where $a_{k} \in\mathbb
{N}, 1 \leq k \leq K$, and $\sum_{k}a_{k} \leq S$. Action $\mathbf{a}$
denotes the deployment of $a_{k}$ servers from the central pool to service
station $k, 1\leq k\leq K$. Should action $\mathbf{a}$ be taken in
state $\mathbf{n}$ then an exponentially distributed amount of time
with rate
%
\begin{equation} \label{e10}
\Lambda= \sum_{k}\{ \lambda_{k} + \mu_{k}( a_{k})
I( n_{k}>0) \}
\end{equation}
will elapse before a change of state. In~(\ref{e10}) $I$ is an
indicator function. The next state of the system
will be $\mathbf{n}+\mathbf{e}_{k}$ with probability $\lambda
_{k}/\Lambda$
and will be $\mathbf{n}-\mathbf{e}_{k}$ with probability $\mu
_{k}(a_{k})I( n_{k}>0) /\Lambda, 1\leq k\leq K$.

A \textit{DSM admissible policy} is given by a map $\mathbf{u}\dvtx\mathbb
{N}^{K}\rightarrow\Xi$, where
%
\begin{equation} \label{e11}
\Xi= \biggl\{ \mathbf{a}; a_{k}\in\mathbb{N},1\leq k\leq K,\mbox{ and
}\sum_{k}a_{k} \leq S\biggr\}
\end{equation}
and is a rule for choosing admissible actions as a function of the current
system state. The cost associated with policy $\mathbf{u}$ is given by
%
\begin{equation} \label{e12}
C(\mathbf{u})=\sum_{k} h_{k}N_{k}(\mathbf{u}),
\end{equation}
where the $h_{k}$ are positive weights (holding cost rates) and
$N_{k}(\mathbf{u})$ is the time average number of customers at station
$k$ under
policy $\mathbf{u}$. The optimization problem of interest is given by
%
\begin{equation} \label{e13}
C^{\mathrm{opt}} = \inf_{\mathbf{u} \in\bar{\mathbf{U}}}C(\mathbf{u}),
\end{equation}
where in~(\ref{e13}) the infimum is over the set $\bar{\mathbf{U}}$ of
DSM admissible policies.

We pause to note that this problem does indeed belong to the class of
dynamic resource allocation problems described in the preceding
section. We
make the choices $c_{k}(a_{k},n_{k})=h_{k}n_{k},
r_{k}(a_{k},n_{k})=a_{k}, L_{k}=S, 1 \leq k \leq K$, with the
transition rates $q_{k}(n_{k}^{\prime}\mid n_{k},a_{k})$ satisfying
\begin{eqnarray*}
q_{k}(n_{k}+1\mid n_{k},a_{k})&=&\lambda_{k}, \\
q_{k}(n_{k}-1\mid n_{k},a_{k})&=&\mu_{k}(a_{k})I(n_{k}>0),
\end{eqnarray*}
for all choices of $k, n_{k}$ and $a_{k}$. They are otherwise zero.
One thing
which is special about this problem is that it is possible to utilize
all of
the resource which is on offer all of the time. It is plainly optimal
to do
so. Hence, in~(\ref{e11}), we can restrict admissible actions to those
which deploy all servers from the pool.

Before proceeding to develop appropriate notions of full
indexability/indi\-ces, we describe assumptions\vadjust{\goodbreak} we shall make about our
service rate functions~$\mu_{k}(\cdot)$. In Assumption~\ref{ass1} we
use the notation
\[
\lceil x \rceil= \min\{ y;y\in\mathbb
{Z}^{+}\mbox{
and }y>x\},\qquad x\in\mathbb{R}^{+}.
\]

\begin{assump} \label{ass1}
There exist functions $\tilde{\mu} _{k}\dvtx\mathbb{R}^{+} \rightarrow
\mathbb{R}^{+}$ which are strictly increasing, twice differentiable and
strictly concave, satisfying
%
\begin{equation} \label{e14}
\tilde{\mu}_{k}(a) = \mu_{k}(a),\qquad a \in[ 0,S ] \cap
\mathbb{N},
\end{equation}
and
%
\begin{equation} \label{e15}
\sum_{k=1}^{K}\lceil\tilde{\mu}_{k}^{-1}(\lambda_{k})
\rceil<S.
\end{equation}
\end{assump}

From~(\ref{e14}) the functions $\tilde{\mu}_{k}, 1 \leq k \leq K$,
are smooth extrapolations of the service rates on the integers in the
range $[ 0,S]$. The properties of these functions reflect
the fact
that, while an increase in the size of the team at a station results in a
higher service rate, the marginal benefit of adding an additional member
diminishes as the team size grows. Requirement~(\ref{e15}) guarantees
the existence of \textit{stable policies} under which all queue
lengths remain finite.
\begin{Remark*}
It is the assumption of strict concavity of the service rate functions at
each station which stimulates an active approach to the distribution of the
pool of servers around the stations and which makes this an interesting
problem. Had we assumed, for example, that the service rates were all convex
in the team size, then \cite{sob82} shows that in an optimal policy the
service pool would always be allocated \textit{en bloc} and we are driven
back to the ``single server'' world of the simple bandit models. This
result is intuitively obvious, as observed by Richard Serfozo to Sobel:
``the fastest rate is also the cheapest.'' Indeed, the resulting
service control problem has a well-known solution in the form of the
so-called $c\mu$-rule. (See \cite{coxsmith61}.)
\end{Remark*}

We are able to develop a Lagrangian relaxation of the problem in (\ref
{e12}) and~(\ref{e13}) as in the preceding section. As in the analysis
of Section~\ref{sec2} up to~(\ref{e5}), such a relaxation yields $K$
optimization problems $P(k,W)$,
one for each station, which here take the form
%
\begin{equation} \label{e16}
C_{k}(W) = \inf_{u_{k}} \{ h_{k}N_{k}(u_{k})+WS_{k}(u_{k})
\},
\end{equation}
where in~(\ref{e16}), the infimum is over the class of DSM
policies $u_{k}\dvtx\mathbb{N} \rightarrow\{ 0,1,\ldots,\break S\}$
which can deploy any number of servers (up to $S$) at station $k$ at
each epoch, $N_{k}(u_{k})$ is the
time average head count and $S_{k}(u_{k})$ the time average number of
servers deployed at $k$ under policy $u_{k}$. The optimization problem in
(\ref{e16}) concerns \textit{station $k$ alone}
and seeks to choose, at each station $k$ decision epoch and in response to
queue length information for station $k$, the number of servers\vadjust{\goodbreak} (from
the $S$
available) to be deployed there. The goal is to make such choices to
minimize costs which are an aggregate of those incurred through customers
waiting $[h_{k}N_{k}(u_{k})]$ and charges imposed for the provision of
service $[WS_{k}(u_{k})]$. Note that Lagrange multiplier $W$ here has an
economic interpretation as the charge imposed per server per unit of time.

We now wish to develop index heuristics for our service allocation
problem by developing station indices of the form described in the
preceding section. These flow from the property of \textit{full
indexability} defined with respect to solutions to the problems
$P(k,W), 1 \leq k\leq K$, and described in Definition~\ref{defn1}.
However, full indexability is a property of individual stations and
hence we now focus on a single station and drop the station identifier
$k$ until further notice. For clarity, the single station problem
$P(W)$ is formulated as an SMDP as follows:
\begin{enumerate}
\item The state of the system at time $t\in\mathbb{R}^{+}$ is $n(t)$,
the number of customers (head count) at the station. New customers
arrive at the station according to a Poisson process of rate $\lambda$.

\item Decision epochs occur at time 0 and whenever there is a change of
state. At each
such epoch an action from the set $A \equiv\{ 0,1,\ldots
,S
\}$
is chosen. Should action $a \in A$ be chosen at time $t$ at which point
$n(t)=n>0$ then costs will be incurred from $t$ at rate $hn+Wa$ and the first
event following $t$ will occur at time $t+X$ where $X \sim\exp[
\lambda+\mu(a)]$. With probabilities $\lambda[\lambda
+\mu
( a) ] ^{-1}$ and $\mu( a) [
\lambda
+\mu
( a) ]^{-1}$ the event will be, respectively, an
arrival or
a service completion.

\item The goal of analysis will be the determination of a stationary policy
to minimize the average cost rate incurred over an infinite horizon.
Trivially, optimal policies offer no service $(a=0)$ when the system is
empty \mbox{$[n(t)=0]$}.
\end{enumerate}

The quest for full indexability is greatly simplified in this case by the
existence of optimal policies for $P(W)$ for which the choice of number
of servers
is increasing in the current head count. We call such policies \textit
{monotone}. This conclusion follows from Theorem 4 in Stidham and Weber
\cite{stidweb89}, which applies to a queueing system with state space
$\mbb{N}$ and Poisson arrivals with an objective which combines a
holding cost which is both increasing in the state and unbounded, with
action costs which are nonnegative and increasing in the resource
level. All of these requirements hold in $P(W)$. Stidham and Weber's
analysis first considers the problem of choosing a policy to minimize
the expected cost incurred in moving the system from a general initial
state to the empty state (their Theorem 2) and then deploys arguments
from renewal theory to demonstrate that such a policy will also
minimize long run average costs (their Section 1.3). We state our
conclusion as Proposition~\ref{prop3}.
\begin{prop}[(Stidham and Weber)] \label{prop3}
There exists a monotone policy which is optimal for $P(W)$.\vadjust{\goodbreak}
\end{prop}

The problem of establishing monotonicity with respect to queue size of
optimal policies for service control problems for queues with Poisson
input is not new. In addition to Stidham and Weber \cite{stidweb89},
see \cite{crab72,doshi78,gall79,geohar01,mit73}. While such
monotonicity is helpful in establishing \textit{full indexability} and in
the subsequent computation of \textit{index functions}, it is not the key
to proving the latter. This is rather the demonstration (to which we
now proceed in Section~\ref{sec31}) that optimal policies for $P(W)$ are
monotone in $W$. Proving this significantly extends the literature on
service control problems for $M/M/1$ queues.

\subsection{Stations are fully indexable}\label{sec31}

In light of Proposition~\ref{prop3} we can recast and simplify the
requirements of
full indexability expressed in Definition~\ref{defn1}. Let $u(W)$ be an
optimal policy
for $P(W)$ which is monotone. It follows that for all choices of $W \in
\mathbb{R}^{+}$ and $0\leq a \leq S-1$,
\[
\Pi\{ u(W),a \} \equiv\{ n \in\mathbb{N}; u(W,n)
\leq a \} = \{ 0,1,\ldots,N(a,W)\}
\]
for some $N(a,W)\in\mathbb{N}\cup\{ \infty\}$. We now
have the following:
\begin{defn}[(Full indexability)] \label{defn3}
The station will be fully indexable if there exists a family of DSM
policies $\{ u(W),W\in\mathbb{R}^{+}\}$ for which (i)
$u(W)$ is monotone and optimal for $P(W)$ $\forall W \in\mathbb
{R}^{+}$ and (ii) the corresponding $N(a,W)$ is increasing in $W,
\forall a \in\{ 0,1,\ldots, S-1\}$.
\end{defn}

To summarize the requirements of Definition~\ref{defn3}, a station will
be fully indexable if the service charge problem $P(W)$ has a monotone
optimal policy for which the number of servers deployed is
\textit{decreasing} in the service charge $W$ for any given head count.
Full indexability enables a \textit{calibration} of the individual
stations as described in Definition~\ref{defn4}.
\begin{defn}[(Station indices)] \label{defn4}
If the station is fully indexable, the corresponding
index function $W\dvtx\{ 0,1,\ldots,S-1\} \times\mathbb{N}
\rightarrow\mathbb{R}^{+}$ is given by
%
\begin{equation} \label{e17}
W(a,n)=\inf\{ W;n\leq N(a,W)\}.
\end{equation}
\end{defn}

In light of Proposition~\ref{prop3} above, Lemma~\ref{lem1} may be
extended as follows in this case:
\begin{lem} \label{lem4}
If the station is fully indexable, the index $W(a,n)$ is \textup{(i)}
decreasing in $a$ for fixed $n$ and \textup{(ii)} increasing in $n$ for
fixed $a$.
\end{lem}

Please note that optimal policies for $P(W)$ will be unchanged if all cost
rates (both holding costs and service charges) are divided by $W>0$
throughout. When we do that, we see that increasing $W$ is equivalent to
decreasing the holding cost rate $h$ in problems for which the service
charge rate is fixed. This being so, we develop the following convenient
reformulation of the definition of\vadjust{\goodbreak} full indexability above: refer to the
problem obtained by setting $W=1$ in the above [namely $P(1)$] as
$Q(h)$ to
emphasize dependence on the holding cost parameter $h$. Hence, $Q(h)$
is the
problem given by
\[
\hat{C}(h)=\inf_{u}\{ hN(u)+S(u)\}.
\]
From Proposition~\ref{prop3} we are able to assert the existence of optimal
policies for $Q(h)$ which are monotone. The following is trivially
equivalent to Definition~\ref{defn3} above.
\begin{defn}[(Full indexability---alternative definition)]
\label{defn5}
The station will be fully indexable if there exists a family of DSM
policies $\{ u(h), h\in\mathbb{R}^{+}\}$ such that, (i)
$u(h)$ is optimal for $Q(h)$  $\forall  h \in\mathbb{R}^{+}$; (ii)
each $u(h)$ is monotone with
\[
\Pi\{ u(h),a\} =\{ 0,1,\ldots,M( a,h)
\},
\]
where $M(a,h)$ is deceasing in $h$  $\forall a \in\{
0,1,\ldots, S-1\}$.
\end{defn}

To summarize, to achieve full indexability, instead of requiring (according
to Definition~\ref{defn3}) that the optimal service level decreases
with the
service charge $W$ (for a fixed value of the holding cost rate $h$), we now
equivalently require it to increase with the holding cost rate $h$ (for
fixed service charge $W=1$). This reformulation of full indexability which
focuses attention on the holding cost element of the objective yields a
more accessible account.

We begin this part of our analysis by noting that it is easy to
establish that \textit{any} optimal policy $u(h)$ for $Q(h)$ must be such
that $\mu\{u(h,n)\} > 0, n \geq1$. It follows that the head count
process is ergodic under its operation. We uniformize station evolution
by rescaling time such that
\[
\lambda+ \mu(S) = 1.
\]
Under this uniformization, the DP
optimality equations for the problem $Q(h)$ are as follows:
%
\begin{eqnarray} \label{e18}\qquad
\lambda v(h,n) & = & hn + \lambda v(h,n+1) \nonumber\\[-8pt]\\[-8pt]
&&{} + \min_{a} \{ a-\mu(a)[
v(h,n)-v(h,n-1)] \} - \gamma(h),\qquad n \geq1 , \nonumber
\end{eqnarray}
where the minimum in~(\ref{e18}) is over the range $0\leq a\leq S$.
Note that in~(\ref{e18}) the quantity $\gamma(h)$ is the minimized
cost rate for $Q(h)$ with $v(h,\cdot)$ the corresponding bias function,
where $v(h,0)=0$. If we write $\hat{C}(h,n,t)$ for the minimum total cost
incurred in $Q(h)$ during $[0,t)$ when $n(0)=n$, then we have $\hat
{C}(h,n,t) \sim t\gamma(h)+v(h,n)$.

Action $a$ is optimal for $Q(h)$ in state $n$ if and only if it
achieves the
minimum in~(\ref{e18}). To proceed further, we write $\Delta v(h,n)
\equiv v(h,n)-v(h,\allowbreak n-1), n \geq1$, and $\Delta v(h,0)=0$. Hence (\ref
{e18}) now becomes
%
\begin{equation} \label{e19}\quad
-\lambda\Delta v(h,n+1) = hn+ \min_{a}\{ a-\mu(a)\Delta
v(h,n)\}
-\gamma(h),\qquad n\geq0.
\end{equation}
We note in passing that it is trivial to deduce from the inductive
specification of $\Delta v(h,\cdot)$ given by the optimality equations,
that the quantities $\{ \Delta v(h,n),n\geq1\}$ are well
defined, including in the event that there are several optimal
policies for~$Q(h)$. The following is an immediate consequence of~(\ref{e19}).
\begin{lem} \label{lem5}
A DSM policy $u$ is optimal for $Q(h)$ if and only if
%
\begin{eqnarray} \label{e20}
&&\Delta v(h,n) \bigl[ \mu\bigl( u(n)+1\bigr) -\mu(u(n))\bigr] \nonumber\\[-8pt]\\[-8pt]
&&\qquad \leq 1
\leq\Delta v(h,n)\bigl[ \mu( u(n))
-\mu
\bigl(u(n)-1\bigr)\bigr],\qquad n \geq1, \nonumber
\end{eqnarray}
where $\mu(S+1)=\mu(S)$ in~(\ref{e20}).
\end{lem}

Please note that if a policy $u$ is such that the inequalities in (\ref
{e20}) are all strict then it is uniquely optimal and so must be monotone
by Proposition~\ref{prop3}. Should the left-hand inequality be
satisfied as an
equation for some $n$ with $u(n)<S$, then both $u(n)$ and $u(n)+1$ are
optimal choices of action in state $n$. To develop the analysis further we
need information regarding the quantities $\Delta v(h,n)$ when viewed as
functions of $h$.
\begin{lem} \label{lem6}
The function $\Delta v(\cdot,n)$ is continuous $\forall n \geq1$.
\end{lem}
\begin{pf}
It is trivial to establish that the average cost rate $\gamma(h)$ is
continuous in $h$. Observe from~(\ref{e19}) that
\[
\Delta v(h,1)=\lambda^{-1} \gamma(h)
\]
and hence $\Delta v(\cdot,1)$ is continuous. From~(\ref{e19}) we
also note that it is straightforward to establish that, if $\Delta
v(\cdot,n)$ is continuous, then so must be $\Delta v(\cdot,n+1),\break
n\geq1$. The result follows by an induction argument.
\end{pf}

Now use $u(h)$ to denote any DSM policy which is optimal for $Q(h)$.
We use $T[u(h),n]$ for the expected time until the system is first emptied
under $u(h)$ given that $n(0)=n$. We also use $C[u(h),n]$ for the
expected cost incurred under $u(h)$ from time 0 when $n(0)=n$ until the
system first empties.
\begin{lem} \label{lem7}
$\forall h>0$,
\[
\Delta v(h,n) \geq\bigl\{ T(u(h),n)-T\bigl(u(h),n-1\bigr)\bigr\} \{
hn-\gamma(h)\} \rightarrow\infty,\qquad n\rightarrow\infty.
\]
\end{lem}
\begin{pf}
A standard argument, based on the fact that the system evolving under $u(h)$
regenerates upon every entry into the empty state, yields the
conclusion that
%
\begin{equation} \label{e21}
v(h,n) =C (u(h),n) - \gamma(h) T(u(h),n),\qquad n \geq1,
\end{equation}
from which we immediately infer that
%
\begin{eqnarray} \label{e22}
\Delta v(h,n) &=& \bigl\{ C(u(h),n)-C\bigl(u(h),n-1\bigr)\bigr\}
\nonumber\\[-8pt]\\[-8pt]
&&{} - \gamma(h) \bigl\{T(u(h),n)-T\bigl(u(h),n-1\bigr)\bigr\},\qquad
n \geq1. \nonumber
\end{eqnarray}

Consider now the system evolving under $u(h)$ from time $0$ when its state
is $n$ until it enters state $n-1$ for the first time. The expected time
taken is plainly $T[u(h),n]-T[u(h),n-1]$ and the holding cost rate incurred
through this period is bounded below by $hn$. If we write the mean
integrated head count divided by $T[u(h),n]-T[u(h),n-1]$ as $\chi[
u(h),n] \geq n$ and
the mean total service cost divided by $T[u(h),n]-T[u(h),n-1]$ as $\psi
[ u(h),n] \geq1$ we infer that
%
\begin{eqnarray} \label{e23}
&&C(u(h),n)-C\bigl(u(h),n-1\bigr) \nonumber\\
&&\qquad= \{h\chi(u(h),n) + \psi(u(h),n)\}
\bigl\{T(u(h),n)-T\bigl(u(h),n-1\bigr)\bigr\}\\
&&\qquad\geq hn \bigl\{ T(u(h),n)-T\bigl(u(h),n-1\bigr)\bigr\},\qquad
n\geq1.\nonumber
\end{eqnarray}
The inequality in the lemma follows immediately from~(\ref{e22}) and
(\ref{e23}). To justify the divergence claim, we simply observe that
an assumed
permanent utilization of the maximum service rate $\mu(S)$ implies that
$\{\mu(S)-\lambda\}^{-1}$ is a
uniform lower bound on $T[u(h),n]-T[u(h),n-1], n \geq1$. The proof
is complete.
\end{pf}

Before proceeding, we observe from~(\ref{e22}) and~(\ref{e23}) and the
definitions of the quantities concerned that we may
write
%
\begin{eqnarray} \label{e24}
\Delta v(h,n) &=& [h\{ \chi(u(h),n)-\alpha(u(h)) \chi
(u(h),1)
\} \nonumber\\
&&\hspace*{1.4pt}{} +\{\psi(u(h),n) - \alpha(u(h)) \psi(u(h),1)\}]
\\
&&{} \times\bigl\{T(u(h),n)-T\bigl(u(h),n-1\bigr)\bigr\},\qquad n\geq1, \nonumber
\end{eqnarray}
where
\[
\alpha(u(h)) = T(u(h),1) [ T(u(h),1) + \lambda^{-1}]^{-1}.
\]
Note that it is straightforward to establish that
%
\begin{equation} \label{e25}
\chi(u(h),n) \geq\chi(u(h),1) > \alpha(u(h)) \chi(u(h),1),\qquad n\geq1.
\end{equation}
The following is an immediate consequence of~(\ref{e20}) and Lemma
\ref{lem7}.
\begin{lem} \label{lem8}
$\forall h>0,  \exists N_{h}< \infty$ such that $u(h,n)=S,
n \geq N_{h}$, for
all choices of $u(h)$.
\end{lem}

We are now in a position to prove full indexability. The key fact to
establish is that $\Delta v(h,n)$ is increasing in $h$ for each $n\geq1$.
Full indexability will then follow trivially from~(\ref{e20}).
\begin{theorem}[(Full indexability)] \label{thm9}
\textup{(i)} The function $\Delta v(\cdot,n)$ is increasing \mbox{$\forall n
\geq1$}; \textup{(ii)} the station is fully indexable.\vadjust{\goodbreak}
\end{theorem}

\begin{pf}
Fix $h_{0}>0$. There are two possibilities. Either there exists a monotone
policy $u(h_{0})$ which is uniquely optimal for $Q(h_{0})$ (case 1) or not
(case $2$).  Under case $1$, invoking the preceding lemma we may
assert the
existence of $N_{h_{0}}<\infty$ such that~(\ref{e20}) is satisfied in the
form
%
\begin{eqnarray} \label{e26}
&&\Delta v(h_{0},n) \bigl[ \mu\bigl( u(h_{0},n)+1\bigr) - \mu
(u(h_{0},n)) \bigr] \nonumber\\
&&\qquad < 1
< \Delta v(h_{0},n)\bigl[ \mu(u(h_{0},n))-\mu
\bigl(u(h_{0},n)-1\bigr)\bigr],\nonumber\\[-8pt]\\[-8pt]
&&\eqntext{1\leq n\leq N_{h_{0}}-1,} \\
&&\hspace*{14.24pt} 1<\Delta v(h_{0},N_{h_{0}})[ \mu(S)-\mu
(S-1)]. \nonumber
\end{eqnarray}
Since $\Delta v(\cdot,n)$ is continuous for $n\geq1$, it must follow
that $\exists \varepsilon>0$ with the property that the inequalities in
(\ref{e26}) are satisfied with $h$ replacing $h_{0}$ for all $h$ in
the range
$h_{0} \leq h < h_{0}+\varepsilon$. We infer from~(\ref{e20}) that monotone
policy $u(h_{0})$ is uniquely optimal for $Q(h), h \in
(h_{0},h_{0}+\varepsilon)$. If we now consider the expression in~(\ref{e24})
with $\alpha, \chi, T$ computed with respect to policy $u(h_{0})$,
it follows
easily that $\Delta v(h,n)$ is increasing and linear in $h$ over the
range $h_{0} \leq h < h_{0}+\varepsilon$.

Now consider case 2. Use $\Upsilon(h_{0})$ to denote the collection of
DSM policies which are optimal for $Q(h_{0})$. From the preceding
lemma and invoking the strict concavity of $\mu(\cdot)$, we infer that
$\Upsilon(h_{0})$ must be finite. Further, the continuity of $\Delta
v(\cdot,n), n\geq1$, together with~(\ref{e20}) implies the existence
of $\delta>0$ such that $Q(h)$ must be optimized by a member of
$\Upsilon(h_{0})$
for $h$ in the range $h_{0}\leq h<h_{0}+\delta$. Suppose that $u\in
\Upsilon(h_{0})$ optimizes $Q(h)$ for some $h\in(h_{0},h_{0}+\delta
)$. It
then follows from~(\ref{e24}) that
%
\begin{eqnarray} \label{e27}
\Delta v(h,n) &=& [h\{ \chi(u,n)-\alpha(u) \chi(u,1)\} +
\{
\psi(u,n) - \alpha(u) \psi(u,1)\}] \nonumber\\[-8pt]\\[-8pt]
&&{} \times\{T(u,n)-T(u,n-1)\},\qquad n\geq1, \nonumber
\end{eqnarray}
where in~(\ref{e27}), $\alpha(u), \chi(u,\cdot), \psi(u,\cdot)$
and $T(u,\cdot)$ denote quantities computed with respect to policy
$u$. Hence
from~(\ref{e27}), it follows that for each $n\geq1$, $\Delta v(\cdot,n)$
lies on one of a finite collection of straight lines with positive gradient
[one for each $u\in\Upsilon(h_{0})$] throughout the range $h_{0} \leq
h < h_{0} + \delta$. However, the continuity of $\Delta v(\cdot,n)$ implies
that it must in fact lie on just one of those lines throughout that range.
It follows that $\Delta v(h,n)$ is increasing linear in $h$ over the
range $h_{0} \leq h < h_{0} + \delta$. We conclude from the above
consideration of cases 1 and 2 that, for each $n\geq
1, \Delta v(\cdot,n)$ is continuous with a positive right gradient at
each $h>0$ and is thus increasing. This concludes the proof of part (i).

For part (ii), we first take the analysis of part (i), case 2, a little
further. Since for the chosen $\delta>0, \Delta v(h,n)$ is strictly
increasing through $[h_{0},h_{0}+\delta)$ for all $n\geq1$, the only
policy which can remain
optimal throughout this range must satisfy conditions of the form (\ref
{e26}). This policy must be maximal (i.e., must assign maximal service
levels) among those policies in $\Upsilon(h_{0})$ and will be uniquely optimal
for $h\in(h_{0},h_{0}+\delta)$ and hence monotone.\vadjust{\goodbreak}

From the above discussion, we can infer the following: fix any
$h_{0}>0$ and
choose the maximal optimal policy for $Q(h_{0})$. This policy is
monotone. Call it $u(h_{0})$. Define $h_{1}$ by
\[
h_{1} = \inf\{h>h_{0};u(h_{0}) \mbox{ is not optimal for } Q(h)\}.
\]
By the above argument $h_{1}>h_{0}$ and $u(h_{0})$ is strictly optimal
for $Q(h), h \in(h_{0},h_{1})$. Further, if $h_{1}<\infty,
u(h_{0})$ is optimal
for $Q(h_{1})$, but not uniquely so. We use $u(h_{1})$ for the maximal
DSM policy which is optimal for $Q(h_{1})$. Policy $u(h_{1})$ is
monotone such that
%
\begin{equation} \label{e27a}
u(h_{1},\cdot) > u(h_{0},\cdot),
\end{equation}
where~(\ref{e27a}) means
\[
u(h_{1},n) \geq u(h_{0},n), \qquad n \geq1,
\]
with strict inequality for at least one $n$. In this way we can develop a
sequence $h_{0} < h_{1} < h_{2} < \cdots< h_{N} < \infty$ and corresponding
monotone policies $u(h_{r}),  0 \leq r\leq N$, such that:
\begin{enumerate}
\item$u(h_{r})$ is optimal for $Q(h),  h\in[ h_{r},h_{r+1}],  0
\leq
r\leq N-1$;

\item$u(h_{r+1},\cdot) > u(h_{r},\cdot),  0 \leq r \leq N-1$;

\item$u(h_{N})$ is optimal for $Q(h),  h \in[ h_{N},\infty)$
and is such that \mbox{$u(h_{N},n)=S,  n\geq1$}.
\end{enumerate}
Since the choice of $h_{0}$ was arbitrary, indexability follows trivially
from 1--3. This completes the proof of part (ii) and of the theorem.
\end{pf}

\subsection{Computation of station indices}\label{sec32}

In the proof of Theorem~\ref{thm9} we constructed an ascending set of
$h$-values,
each of which signaled a change of optimal policy for~$Q(h)$. In this
construction the initial $h_{0}$ was arbitrary. In our discussion of index
computation, we shall continue initially to operate in $h$-space [i.e., to
consider solutions to the optimization problems $Q(h)$], but will
construct a \textit{descending} set of $h$-values, labeled
$j_{1},j_{2},\ldots$ each of which will also signal a change of optimal
policy. We do this because such a set is straightforward to initialize, with
$j_{1}$ the supremum of those $h$ for which the policy [hereafter
labeled $u(j_{0})$] which applies the maximal number of servers $S$
whenever the
queue is nonempty is \textit{not} optimal for $Q(h)$. Because of our
ability to restrict to monotone policies, it is clear that both $u(j_{0})$
and the policy $u(j_{1})$ (which applies $S-1$ servers when the queue length
is~$1$, but which otherwise applies $S$ servers) are optimal for
$Q(j_{1})$. By direct calculation of the average cost rates for these
policies it is straightforward to verify that
\[
j_{1}=\{ \mu(S)-\lambda\} \biggl\{ \frac{1}{\mu(S)-\mu
(S-1)}-\frac{S}{\mu(S)}\biggr\}.
\]
We now give an algorithm for producing the sequence $\{
j_{m},m\geq
1\}$ and the monotone policies $\{ u(j_{m}),m\geq0
\}$
such that $u(j_{m})$ is strictly optimal for $Q(h)$ in the range
$j_{m+1}<h<j_{m}$. Note that we take $j_{0}=\infty$.\vadjust{\goodbreak} In the algorithm we
utilize the characterization of optimal policies for $Q(h)$ given in
Lemma~\ref{lem5}
together with the formula for $\Delta v(h,n)$ given following the proof
of Lemma~\ref{lem7}.

\subsubsection*{Algorithm for index computation}

\mbox{}

\textit{Step} 0. Let $m=1$. The positive real $j_{1}$ and
the policy $u(j_{1})$ are as
above. The positive integer $N_{1}$ is given by
\[
N_{1}=\min\{ n;u(j_{1},n)=S\} =2.
\]

\textit{Step} 1. The positive real $j_{m}$, the policy
$u(j_{m})$ and the positive
integer $N_{m}$ given by
\[
N_{m}=\min\{ n;u(j_{m},n)=S\}
\]
are specified. Determine $(A_{n}^{m},B_{n}^{m};1\leq n\leq N_{m})$
given by
\[
A_{n}^{m} = \{ \chi(u(j_{m}),n) - \alpha(u(j_{m})) \chi
(u(j_{m}),1) \} \bigl\{ T(u(j_{m}),n)-T\bigl(u(j_{m}),n-1\bigr)\bigr\}
\]
and
\[
B_{n}^{m} = \{ \psi(u(j_{m}),n) - \alpha(u(j_{m})) \psi
(u(j_{m}),1)\} \bigl\{ T(u(j_{m}),n)-T\bigl(u(j_{m}),n-1\bigr)\bigr\}.
\]

\textit{Step} 2. Let $j_{m+1}$ be the maximal $h$ satisfying
\[
\{ A_{n}^{m}h + B_{n}^{m}\} \bigl\{\mu(u(j_{m},n)) - \mu
\bigl(u(j_{m},n)-1\bigr)\bigr\} = 1
\]
for some $n$ in the range $1\leq n\leq N_{m}$. Let $n_{m}$ be an $n$-value
achieving the equality.

\textit{Step} 3. Define the policy $u(j_{m+1})$ by
\[
u(j_{m+1},n)=u(j_{m},n)-I(n=n_{m}),\qquad n\geq0,
\]
where $I$ is an indicator. Determine $N_{m+1}$ and the
$(A_{n}^{m+1},B_{n}^{m+1};1\leq n\leq N_{m+1})$ as in Step 1.

\textit{Step} 4. If $j_{m+1}\leq0$, stop. Otherwise
return to Step 2.

It is now straightforward to recover the station indices (as given in
Definition~\ref{defn2}) from the quantities calculated by the above
algorithm. Note,
as previously, that optimal policies for $P(W)$ and $Q(h/W)$
coincide. In order to compute the station index $W(a,n)$, determine
from the
above algorithm the value $j_{m}$ satisfying
\[
u(j_{m},n) = a+1 \quad\mbox{and}\quad u(j_{m+1},n) = a.
\]
We then infer that
\[
W(a,n)=\frac{h}{j_{m+1}}.
\]

\subsection{Numerical study}\label{sec33}

Extensive numerical investigations have been conducted on the
performance of the greedy index heuristic as a policy for the queueing
control problems described above. We shall now present some of our
results as Examples~\ref{example1} and~\ref{example2}.\vadjust{\goodbreak}
\begin{example}\label{example1}
All Example~\ref{example1} problems concern the dynamic allocation of a pool of
twenty-five servers ($S=25$) to two service stations ($K=2$). Service
rate functions have the form
%
\begin{equation} \label{e34}
\mu_k(a) = a( a + \nu_k)^{-1} \mu_k,\qquad k =1,2.
\end{equation}
In all, 4950 problems were generated at random, consisting of 99 sets
of 50 problems. For each problem the parameters $\lambda_1, \lambda
_2, \mu_1, \mu_2, \nu_1, \nu_2$ were chosen by sampling
independently from uniform distributions. Full details may be found at
\url{http://www.lums.lancs.ac.uk/files/onlinesup.pdf}.

\begin{table}
\caption{Choices of the parameters $\lambda_1, \lambda_2, \mu_1$ and
$\mu_2$ $(G, J)$ and $\nu_1, \nu_2$ $(7)$ and $\eta_1, \eta_2$ $(14)$
which give challenging problem sets for Examples \protect\ref{example1} and
\protect\ref{example2}} \label{tab1}
\begin{tabular*}{\tablewidth}{@{\extracolsep{\fill}}lccc@{}}
\hline
\multicolumn{1}{@{}l}{\textit{\textbf{G}}} & \multicolumn{1}{c}{\textit{\textbf{J}}} &
\multicolumn{1}{c}{\textbf{7}} & \multicolumn{1}{c@{}}{\textbf{14}} \\
\hline
$\lambda_1 \in[0.8,1.1)$ & $\lambda_1 \in[0.8,1.1)$ & $\nu_1 \in
[5.0,10.0)$ & $\eta_1 \in[0.07,0.125)$ \\
$\lambda_2 \in[1.6,2.2)$ & $\lambda_2 \in[1.6,2.2)$ & $\nu_2 \in
[0.5,2.0)$ & $\eta_2 \in[0.2,0.3)$ \\
$\mu_1 \in[1.5,1.8)$ & $\mu_1 \in[1.5,1.8)$ &\\
$\mu_2 \in[3.0,3.6)$ & $\mu_2 \in[4.4,5.0)$ &\\
\hline
\end{tabular*}
\end{table}

\begin{table}[b]
\caption{The percentage cost rate excess over optimum of \textup{(i)} the
greedy index heuristic for~all~$4950$ Example~\protect\ref{example1} problems, \textup{(ii)} for
problem sets $G7$ and $J7$ and \textup{(iii)}~for~the~best~static~allocation
policy} \label{tab3}
\begin{tabular*}{\tablewidth}{@{\extracolsep{\fill}}lcccd{2.4}@{}}
\hline
& \textbf{Overall} & \textbf{\textit{G}7} & \textbf{\textit{J}7} & \multicolumn{1}{c@{}}{\textbf{Static}} \\
\hline
\textit{MIN} & 0.0000 & 0.0416 & 0.0263 & 1.7837\\
\textit{LQ} & 0.0001 & 0.0745 & 0.0558 & 5.6978\\
\textit{MED} & 0.0021 & 0.0964 & 0.1021 & 8.1880\\
\textit{UQ} & 0.0186 & 0.1670 & 0.1433 & 10.9678\\
\textit{MAX} & 0.2910 & 0.2910 & 0.2422 & 22.1868\\
[4pt]
\textit{N} & 4950 & 50 & 50 & \multicolumn{1}{c@{}}{4950}\\
\hline
\end{tabular*}
\end{table}

For each of the 4950 problems generated, indices were developed using
the algorithm given in Section~\ref{sec32}. Time average holding cost rates
for the greedy index heuristic and an optimal policy were computed
using DP value iteration and the \textit{percentage cost rate excess} of
the index heuristic over the optimum was recorded. Order statistics
(minimum, lower quartile, median, upper quartile, maximum) of the
percentage cost rate excess over optimum of the index heuristic are
given in Table~\ref{tab3} for the 4950 problems overall, together with
those for two of the problem sets ($G7,J7$) for which the heuristic
performed \textit{relatively} less well. For ease of reference, Table
\ref{tab1} gives details of the uniform distributions used to generate
these challenging problem sets. Additionally, in Table~\ref{tab3} under
the head ``Static'' are recorded the order statistics for the
percentage cost rate excess over optimum for the best static policy
which makes a fixed allocation of servers to stations for all time.
These latter values give an indication of the potential value of
designing a dynamic policy for these resource allocation problems.

The greedy index heuristic performs outstandingly well with a worst
case suboptimality of 0.2910\% for one of the problems generated as
part of the problem set $G7$. Inspection of the results for $G7$ and $J7$
show that the performance of the index policy is excellent even in
problems for which the stochastic dynamics of the two stations are very
different. Perusal of the results for individual problems suggests that
the benefits of designing a dynamic policy tend to be greatest when the
greedy index heuristic performs \textit{relatively} less well. For one
particular problem not recorded in Table~\ref{tab3} for which the
greedy index heuristic had a cost rate excess over optimal of 0.8801\%
that of the best static policy was 48.9693\%.
\end{example}
\begin{example}\label{example2}
All Example~\ref{example2} problems concern the dynamic allocation of a pool of
twenty-five servers ($S=25$) to two service stations ($K=2$). Service
rate functions have the form
\[
\mu_k(a) = \bigl(1- \exp(-a\eta_k)\bigr)\mu_k,\qquad k=1,2.
\]
Other details are similar to those of Example~\ref{example1}. Again, 4950 problems
were generated at random in 99 sets of 50. For each problem the
parameters $\lambda_1, \lambda_2, \mu_1, \mu_2$, $\eta_1, \eta_2$ were
chosen by sampling independently from uniform distributions. While
Table~\ref{tab1} gives details of the distributions used for some of
the more challenging problems ($G14, J14$), full details may be found at
\url{http://www.lums.lancs.ac.uk/files/onlinesup.pdf}.

\begin{table}[b]
\caption{The percentage cost rate excess over optimum of \textup{(i)} the
greedy index heuristic for~all~$4950$ Example~\protect\ref{example2} problems, \textup{(ii)} for
problem sets $G14$ and $J14$ and \textup{(iii)}~for~the~best~static~allocation
policy} \label{tab5}
\begin{tabular*}{\tablewidth}{@{\extracolsep{\fill}}lcccd{2.4}@{}}
\hline
& \textbf{Overall} & \textbf{\textit{G}14} & \textbf{\textit{J}14} & \multicolumn{1}{c@{}}{\textbf{Static}}\\
\hline
\textit{MIN} & 0.0000 & 0.0803 & 0.0279 & 2.2079\\
\textit{LQ} & 0.0024 & 0.1473 & 0.1100 & 7.0473\\
\textit{MED} & 0.0087 & 0.2164 & 0.1495 & 10.2092\\
\textit{UQ} & 0.0372 & 0.4289 & 0.2509 & 14.4034\\
\textit{MAX} & 0.8469 & 0.8469 & 0.5905 & 26.5599\\
[4pt]
\textit{N} & 4950 & 50 & 50 & \multicolumn{1}{c@{}}{4950}\\
\hline
\end{tabular*}
\end{table}

For each of the 4950 problems generated, the percentage cost rate
excess of the greedy index heuristic over the optimum was computed. The
overall results are presented in Table~\ref{tab5} along with those for
problem sets $G14$ and $J14$ and for the best static policy. Similar
comments apply as for Example~\ref{example1}.
\end{example}

\section{Spinning plates: Optimal investment in a collection of reward
generating assets}\label{sec4}

As a further illustration of the applicability of the methodology of
Section~\ref{sec2}, we now give a brief account of a setup in which a collection
of $K$ reward generating assets is maintained using a divisible
investment resource. Each asset $k$ evolves on its (finite) state space
$\{0,1,\ldots,A_k\}$ with higher-valued states being those in which the
reward performance of the asset is stronger. In the absence of
investment, assets tend to deteriorate toward lower-valued states.
Positive investment \textit{both} arrests the asset's tendency to
deteriorate and enhances asset performance by enabling upward movement
of the asset state. Investment decisions will often need to strike a
balance between maintaining the performance of highly performing assets
and improving the performance of poorly performing ones. Our model
class represents a significant generalization of the \textit{spinning
plates} model of asset management discussed by Glazebrook, Kirkbride and Ruiz-Hernandez
\cite{gla06} to the case of a divisible resource.

Formally, the \textit{system state} at time $t$ is $\mb{n}(t) = \{
n_1(t),n_2(t), \ldots, n_K(t)\}$, where $n_k(t)$ is the state of asset
$k$ at $t$. The system state is observed continuously with \textit
{decision epochs} at time zero and at subsequent times at which the
system state changes. An \textit{admissible action} is a vector $\mb{a} =
(a_1, a_2, \ldots, a_K)$, with $a_k$ identified with the rate at which
investment resource is applied to asset $k$, where $a_k \in\mbb{N},
1 \leq k \leq K$, and $\sum_k a_k \leq R$. Note that $R$ is the rate at
which investment resource is available to the system.

Functions $\lambda_k\dvtx \{0,1,\ldots,R\} \times\{0,1,\ldots,A_k-1\}
\rightarrow
\mbb{R}^+$ and $\mu_k\dvtx \{0,1,\ldots,R\}$ $\times\{1,2,\ldots,A_k\}
\rightarrow
\mbb{R}^+$ are used in the specification of the transition law of asset
$k$ as follows:
%
\begin{eqnarray} \label{n35}
q_k(n_k+1\mid n_k,a_k) &=& \lambda_k(a_k,n_k) I(n_k < A_k)
\nonumber\\[-8pt]\\[-8pt]
&&\eqntext{\mbox{(Investment enhances asset performance)}}
\end{eqnarray}
and
%
\begin{eqnarray} \label{n36}
q_k(n_k-1\mid n_k,a_k) &=& \mu_k(a_k,n_k) I(n_k > 0)
\nonumber\\[-8pt]\\[-8pt]
&&\eqntext{\mbox{(Investment arrests asset deterioration)}.}
\end{eqnarray}
All other transition rates for asset $k$ are zero. We shall assume that
$\lambda_k(\cdot,n_k)$ is strictly increasing and strictly concave
$\forall n_k \in\{0,1,\ldots, A_k-1\}$ and that $\mu_k(\cdot,n_k)$
is strictly decreasing and strictly convex $\forall n_k \in\{
1,2,\ldots, A_k\}$. These conditions describe laws of diminishing
returns as the level of investment to an asset increases, regardless of
its state. It would be natural in many application contexts to further
assume that each $\lambda_k (a_k, \cdot)$ is decreasing and each $\mu_k
(a_k, \cdot)$ is increasing $\forall a_k \in\{0,1,\ldots,R\}$,
namely, that when an asset is in a higher-valued state, improvements
take longer to achieve but asset deterioration occurs more rapidly. Our
theoretical results do \textit{not} require these latter conditions to
hold, though they will feature in the problems analyzed in our\vadjust{\goodbreak}
numerical study. Finally, in state $\mb{n}$, each asset $k$ earns
returns at rate $d_k(n_k)$, where $d_k\dvtx \{0,1,\ldots,A_k\} \rightarrow
\mbb
{R}^+$ is increasing. The dynamic resource allocation problem of
interest is expressed as
%
\begin{equation} \label{n37}
D^{\mathrm{opt}} = \sup_{\mb{u} \in\bar{\mb{U}}} \sum_k D_k(\mb{u}),
\end{equation}
while in~(\ref{n37}), $\bar{\mb{U}}$ is the set of DSM and admissible
policies and $D_k(\mb{u})$ is the reward rate earned by asset $k$ under
policy $\mb{u}$.

\subsection{Assets are fully indexable}\label{sec41}

Following a version of the development of Section~\ref{sec2} which focuses on
reward maximization instead of cost minimization, we develop a
Lagrangian relaxation of~(\ref{n37}) which yields $K$ single asset
problems $P(k,W), 1 \leq k \leq K$, of the form
%
\begin{equation} \label{n38}
\sup_{u_k} \{D_k(u_k)-WR_k(u_k)\}.
\end{equation}
In~(\ref{n38}), the supremum is over the class of DSM policies $u_k
\dvtx
\{
0,1,\ldots, A_k\} \rightarrow\{0,1,\ldots,R\}$ which can apply any resource
level at asset $k$. Further, $D_k(u_k)$ is the asset $k$ return rate
under policy $u_k$, while $R_k(u_k)$ is the rate of resource consumed.
Full indexability of project $k$ requires the existence of optimal
policies for~(\ref{n38}) which, in every state, apply a resource rate
to the asset which is decreasing in the resource charge $W$. In
discussing full indexability, we now drop the asset subscript $k$ and
use $P(W)$ for the single asset problem
%
\begin{equation} \label{n39}
\sup_u \{ D(u) - W R(u)\}.
\end{equation}
Following the approach of Section~\ref{sec31} we introduce the problem
$Q(h)$, defined by
%
\begin{equation} \label{n40}
\sup_u \{h D(u) - R(u) \}
\end{equation}
and argue that full indexability will be established by the existence
of optimal policies for~(\ref{n40}) which, in every state, choose
resource levels which are increasing in the reward multiplier $h$.

In order to develop the DP optimality equations for $Q(h)$ we
uniformize asset evolution by rescaling time such that
%
\begin{equation} \label{n41}
\max_{0 \leq n \leq A} \{ \lambda(R,n) + \mu(0,n)\} = 1.
\end{equation}
Under the rescaling in~(\ref{n41}), we use $\gamma(h)$ for the maximal
reward rate for $Q(h)$ and $v(h,\cdot)$ for the corresponding bias
function. The optimality equations may be written
%
\begin{eqnarray} \label{n42}
0 &=& -\gamma(h) + hd(n) \nonumber\\
&&{}+ \max_a [ -a + \lambda(a,n) \Delta v(h,n+1)
I(n<A)\nonumber\\[-8pt]\\[-8pt]
&&\hspace*{69.5pt}{} - \mu(a,n)
\Delta v(h,n) I(n>0) ], \nonumber\\
&&\eqntext{0 \leq n \leq A.}
\end{eqnarray}
In~(\ref{n42}), we take $\Delta v(h,n) \equiv v(h,n) - v(h,n-1),
1 \leq n \leq A$, and the maximization is over $0 \leq a \leq R$. Lemma
\ref{lem11} uses~(\ref{n42}) to give a characterization of optimal
policies for $Q(h)$.
\begin{lem} \label{lem11}
A DSM policy $u$ is optimal for $Q(h)$ if and only if
%
\begin{eqnarray} \label{n43}
&&\Delta v(h,n+1) I(n<A) \bigl[ \lambda\bigl(u(n)+1,n\bigr)-\lambda(u(n),n) \bigr]
\nonumber\\
&&\quad{} + \Delta v(h,n) I(n>0) \bigl[\mu(u(n),n)-\mu\bigl(u(n)+1,n\bigr)\bigr] \nonumber\\
&&\qquad \leq1 \leq\Delta v(h,n+1) I(n<A) \bigl[\lambda(u(n),n) - \lambda
\bigl(u(n)-1,n\bigr)\bigr] \\
&&\qquad\quad\hspace*{16.4pt}{} + \Delta v(h,n) I(n>0) \bigl[\mu\bigl(u(n)-1,n\bigr)-\mu(u(n),n)\bigr],\nonumber\\
&&\eqntext{0 \leq n \leq A,}
\end{eqnarray}
where in~(\ref{n43}) we take $\lambda(R+1,\cdot) \equiv\lambda (R,\cdot
), \lambda(-1,\cdot) \equiv- \infty, \mu(R+1,\cdot) \equiv \mu
(R,\cdot), \mu(-1,\cdot) \equiv\infty$.
\end{lem}
\begin{Remark*}
One important point of difference between our generalized spinning
plates model and the queueing models of Section~\ref{sec3} is that the existence
of optimal policies for $Q(h)$ which are monotone in state is no longer
guaranteed, even for assets for which the transition rates are state
monotone for any fixed resource level. Indeed, counter-examples are
easy to find. The following asset appeared in the very first of 2000
randomly generated problems contributing to Table~\ref{tabs5}, which
appears later in Section~\ref{sec42} as part of an extensive numerical
investigation into the performance of the greedy index heuristic.

\begin{table}[b]
\caption{Values of optimal actions (resource levels) for $Q(h)$ for
seven $h$-values and all~states~$0$~(leftmost entry) to $10$
(rightmost entry)} \label{tabs4}
\begin{tabular*}{\tablewidth}{@{\extracolsep{\fill}}lccccccccccc@{}}
\hline
3&4&4&4&3&3&2&2&2&1&0&$h=7.37491$\\
2&4&4&4&3&3&2&2&2&1&0&$h=7.07632$\\
2&4&4&3&3&3&2&2&2&1&0&$h=5.32243$\\
2&4&3&3&3&3&2&2&2&1&0&$h=5.21572$\\
2&3&3&3&3&3&2&2&2&1&0&$h=4.98366$\\
2&3&3&3&3&2&2&2&2&1&0&$h=3.84063$\\
1&3&3&3&3&2&2&2&2&1&0&$h=3.48775$\\
\hline
\end{tabular*}
\end{table}

We make the following asset choices: $R=5, A=10$
\[
\lambda(a,n) = a(a + \phi)^{-1},\qquad 0 \leq a \leq5, 0 \leq n \leq9,
\]
and
\[
\mu(a,n) = \phi(a+\phi)^{-1} \eta,\qquad 0 \leq a \leq5, 1 \leq n \leq10,
\]
where $\phi= 1.30738$ and $\eta= 1.16393$. Further, the return for
the asset is given by $d(n) = n(n+1)^{-1}$. In
Table~\ref{tabs4},\vadjust{\goodbreak}
find values of $u(h,n), 0 \leq n \leq10$, for seven distinct values
of $h$, where $u(h,\cdot)$ is an optimal policy for $Q(h)$. Note that
for the six open $h$-intervals whose endpoints are the successive $h$
values given in Table~\ref{tabs4}, the policy which sits alongside the
value of $h$ which is the lower endpoint is uniquely optimal throughout
the interval. At no value of $h$ in the range $(3.48775,7.37491)$ is
there an optimal policy for $Q(h)$ which is monotone in state. Please
note that the values in Table~\ref{tabs4} are consistent with the
asset's \textit{full indexability} in that optimal actions for any given
state are everywhere increasing in $h$ over the range considered.
\end{Remark*}

We now consider the state process $\{n(t), t \geq0\}$ of a single
asset evolving under some fixed DSM policy $u$ for $Q(h)$. We shall
write $\gamma(u,h)$ for the reward rate earned under policy $u$ and
$v(u,h,\cdot)$ for the corresponding bias function. Recall our earlier
notational choices: if $u(h)$ is optimal for $Q(h)$ then $\gamma
(u(h),h) \equiv\gamma(h)$ and $v(u(h),h,\cdot) \equiv v(h,\cdot)$.

Suppose now that $n(0) = n \in[1,A]$. We define the stopping times
\mbox{$\tau(u,m\mid n)$} by
\[
\tau(u,m\mid n) = \inf\{t>0; n(t) = m\},\qquad 0 \leq m < n \leq A,
\]
to be the first time after time 0 at which the asset state enters $m$
when policy $u$ is applied throughout. We use
%
\begin{eqnarray} \label{n46}
D(u,h,n) &=& h \mbb{E} \biggl[ \int_{0}^{\tau(u,0\mid n)} d\{n(t)\}\,dt\biggr]
- \mbb{E} \biggl[ \int_0^{\tau(u,0\mid n)} u\{n(t)\} \,dt \biggr] \\
\label{n47}
&\equiv& h \chi(u,n) - \psi(u,n),\qquad 1 \leq n \leq A,
\end{eqnarray}
for the expected reward (net of resource charges) earned by the asset
evolving under policy $u$ during $[0,\tau(u,0\mid n))$ and
%
\begin{equation} \label{n48}
T(u,n) = \mbb{E}\{\tau(u,0\mid n)\},\qquad 1 \leq n \leq A.
\end{equation}
As in the proof of Lemma~\ref{lem7} we can use standard renewal
arguments to infer that
%
\begin{equation} \label{n49}
v(u,h,n) = D(u,h,n) - \gamma(u,h) T(u,n),\qquad 1 \leq n \leq A,
\end{equation}
and hence that
%
\begin{eqnarray} \label{n50}
\Delta v(u,h,n) &=& \{D(u,h,n)-D(u,h,n-1)\} \nonumber\\[-8pt]\\[-8pt]
&&{} - \gamma(u,h)\{T(u,n)-T(u,n-1)\},\qquad 1 \leq n \leq A.\nonumber
\end{eqnarray}
We now observe that taking $n=1$ in~(\ref{n46})--(\ref{n48}) yields
%
\begin{eqnarray} \label{n51}
\gamma(u,h) &=& [ h \chi(u,1) - \psi(u,1) + \{hd(0)-u(0)\}\{
\lambda
(u(0),0)\}^{-1} ] \nonumber\\[-8pt]\\[-8pt]
&&{} \times[ T(u,1) + \{ \lambda(u(0),0)\}^{-1}
]^{-1}.\nonumber
\end{eqnarray}
Using~(\ref{n51}) in~(\ref{n50}) we observe that, for any fixed $u,
n$ where $1 \leq n \leq A, \Delta v(u,\break h,n)$ is affine in $h$ with
$h$-gradient proportional to
\begin{eqnarray*}
&& \frac{\chi(u,n)-\chi(u,n-1)}{T(u,n)-T(u,n-1)} - \frac{\chi
(u,1)+d(0)\{\lambda(u(0),0)\}^{-1}}{T(u,1)+ \{\lambda(u(0),0)\}^{-1}}
\\[-2pt]
&&\qquad= \frac{\mbb{E}[\int_0^{\tau(u,n-1\mid n)} d\{n(t)\} \,dt
]}{\mbb
{E}\{\tau(u,n-1\mid n)\}} - \frac{\mbb{E}[\int_0^{\tau(u,0\mid 1)}
d\{
n(t)\}\,dt]+d(0)\{\lambda(u(0),0)\}^{-1}}{\mbb{E}\{\tau
(u,0\mid 1)\} +
\{\lambda(u(0),0)\}^{-1}} \\[-2pt]
&&\qquad\geq\frac{\mbb{E}[\int_{0}^{\tau(u,n-1\mid n)}d\{n(t)\}\,dt
]}{\mbb{E}\{\tau(u,n-1\mid n)\}} - \frac{\mbb{E} [ \int
_{0}^{\tau
(u,0\mid 1)} d\{n(t)\}\,dt]}{\mbb{E}\{\tau(u,0\mid 1)\}},\qquad 1 \leq n \leq A,
\end{eqnarray*}
which is easily seen to be positive since the return rate $d(\cdot)$ is
increasing in the state. We infer that $\Delta v(u,\cdot,n)$ is
increasing for any fixed $u, n$ where $1 \leq n \leq A$. It must,
therefore, follow that $\Delta v(\cdot,n)$ is increasing over any
$h$-interval for which there exists some fixed policy $u(h)$ which is
strictly optimal for $Q(h)$.

We can now deploy arguments along the lines of those in the proof of
Theorem~\ref{thm9} to infer Theorem~\ref{thm12}(i). Please note that
Theorem~\ref{thm12}(ii) follows straightforward from Theorem \ref
{thm12}(i) together with Lemma~\ref{lem11} and the conditions on the
transition rates enunciated after~(\ref{n36}). This result generalizes
Theorem 1 of Glazebrook, Kirkbride and Ruiz-Hernandez \cite{gla06} to the divisible
resource case.\vspace*{-2pt}
\begin{theorem}[(Full indexability)] \label{thm12}
\textup{(i)} The functions $\Delta v(\cdot,n)$ are increasing $\forall n,
1 \leq n \leq A$; \textup{(ii)} the asset is fully indexable.\vspace*{-2pt}
\end{theorem}

We apply an algorithm similar to that in Section~\ref{sec32} to infer the
sequence of optimal policies as $h$ decreases from some large value for
which the optimal policy uses maximal resource $R$ in every state below
$A$. Indices are now \textit{not} in general monotone in state.\vspace*{-2pt}

\subsection{Numerical study}\label{sec42}

We proceed to assess the quality of the greedy index heuristic through
a study of 14,000 randomly generated two asset problems $(K=2)$ in
which resource is available to the assets in integer amounts up to a
maximum of 5 or 10 ($R=5$ or 10). All assets studied evolve over the
state space $\{0,1, \ldots, 10\}$ while the transition rates for each
asset $k$ are assumed to be multiplicatively separable such that
%
\begin{equation} \label{s50}
\lambda_k(a_k,n_k) = a_k(a_k+\phi_k)^{-1} \xi_k(n_k),\qquad 0 \leq a_k
\leq
R, 0 \leq n_k \leq9,
\end{equation}
and
%
\begin{equation} \label{s51}
\mu_k(a_k,n_k) = \phi_k(a_k+\phi_k)^{-1}\eta_k(n_k),\qquad
0 \leq a_k \leq
R,  1 \leq n_k \leq10,
\end{equation}
with $\phi_k$ a positive constant. In all 14,000 problems the $\phi_k$
will be obtained by sampling from the uniform distribution on
$[0.75,5.00]$. The assets are assumed always to have a common return
function, denoted $d\dvtx\{0,1,\ldots,10\} \rightarrow\mbb{R}^+$, which
is increasing.\vadjust{\goodbreak}

\begin{table}[b]
\def\arraystretch{0.92}
\vspace*{-3pt}
\tablewidth=230pt
\caption{The percentage return rate below optimum of \textup{(i)} the greedy
index heuristic, \textup{(ii)} the best static allocation policy and \textup{(iii)}
a~myopic policy for $2000$ problems with state independent transition
rates. See text for details} \label{tabs5}
\begin{tabular*}{\tablewidth}{@{\extracolsep{\fill}}lcd{2.4}d{2.4}@{}}
\hline
& \textbf{Index} & \multicolumn{1}{c}{\textbf{Static}} & \multicolumn{1}{c@{}}{\textbf{Myopic}} \\
\hline
\textit{MIN} & 0.0000 & 0.0719 & 0.0027 \\
\textit{LQ} & 0.1482 & 3.7812 & 4.7774 \\
\textit{MED} & 0.6752 & 6.1724 & 16.7270 \\
\textit{UQ} & 1.0751 & 7.4822 & 26.5042 \\
\textit{MAX} & 1.9082 & 13.6966 & 39.3193 \\
[4pt]
\textit{N} & 2000 & \multicolumn{1}{c}{2000} & \multicolumn{1}{c@{}}{2000}\\
\hline
\end{tabular*}
\end{table}

In all problems we compare the performance of three heuristic policies
for resource allocation. These are the greedy index policy (Index), the
optimal static policy (Static) and a myopic policy (Myopic) which in
every system state $\mb{n}=(n_1,n_2)$ chooses an action $\mb{a} =
(a_1,a_2)$ to maximize the rate at which the return rate from the
assets increases, namely,
\begin{eqnarray*}
&&\max_{\mb{a}} \sum_{k=1}^2 [\lambda_k(a_k,n_k) I(n_k<10) \{
d(n_k+1)-d(n_k)\} \\[-2pt]
&&\hspace*{14.4pt}\qquad{} + \mu_k(a_k,n_k) I(n_k>0) \{d(n_k-1)-d(n_k)\}].
\end{eqnarray*}
For each problem instance, the return rate achieved under each
heuristic is compared to optimum and reported as a percentage
suboptimality. All computations utilize DP value iteration. The
problems are generated in seven groups with 2000 problems in each
group. For each group of problems and each heuristic the 2000
percentage suboptimalities are summarized using order statistics, as
was done in Section~\ref{sec33}. The results are presented in Tables~\ref
{tabs5}--\ref{tabs8}. The problem details now follow.

\begin{table}
\tablewidth=250pt
\caption{The percentage return rate below optimum of \textup{(i)} the greedy
index heuristic, \textup{(ii)} the best static allocation policy and
\textup{(iii)}~a~myopic policy for $2000$ problems with state dependent transition
rates. See text for details} \label{tabs6}
\begin{tabular*}{\tablewidth}{@{\extracolsep{\fill}}lccd{2.4}@{}}
\hline
& \textbf{Index} & \textbf{Static} & \multicolumn{1}{c@{}}{\textbf{Myopic}} \\
\hline
\textit{MIN} & 0.0000 & 0.0305 & 2.2993 \\
\textit{LQ} & 0.0000 & 0.0695 & 6.7075 \\
\textit{MED} & 0.0001 & 0.1179 & 13.0721 \\
\textit{UQ} & 0.0008 & 0.1888 & 17.9062 \\
\textit{MAX} & 0.9685 & 1.0340 & 23.0439 \\
[4pt]
\textit{N} & 2000 & 2000 & \multicolumn{1}{c@{}}{2000}\\
\hline
\end{tabular*}
\end{table}

\begin{table}[b]
\tablewidth=250pt
\caption{The percentage return rate below optimum of \textup{(i)} the greedy
index heuristic, \textup{(ii)} the best static allocation policy and
\textup{(iii)}~a~myopic policy for $2000$ problems with state independent transition
rates. See text for details} \label{tabs7}
\begin{tabular*}{\tablewidth}{@{\extracolsep{\fill}}lccd{2.4}@{}}
\hline
& \textbf{Index} & \textbf{Static} & \multicolumn{1}{c@{}}{\textbf{Myopic}} \\
\hline
\textit{MIN} & 0.0000 & 0.1830 & 1.2736 \\
\textit{LQ} & 0.0000 & 0.3275 & 1.7252 \\
\textit{MED} & 0.0001 & 0.3817 & 1.9311 \\
\textit{UQ} & 0.0012 & 0.4652 & 2.5708 \\
\textit{MAX} & 0.0095 & 0.7310 & 16.1912 \\
[4pt]
\textit{N} & 2000 & 2000 & \multicolumn{1}{c@{}}{2000}\\
\hline
\end{tabular*}
\end{table}

\begin{table}
\caption{The percentage return rate below optimum of \textup{(i)} the greedy
index heuristic, \textup{(ii)}~the~best~static allocation policy and \textup{(iii)} a
myopic policy for $8000$ problems with~state~dependent transition
rates. See text for details} \label{tabs8}
\begin{tabular*}{\tablewidth}{@{\extracolsep{\fill}}ld{2.4}d{2.4}d{2.4}d{2.4}d{2.4}d{2.4}@{}}
\hline
& \multicolumn{1}{c}{\textbf{Index}} & \multicolumn{1}{c}{\textbf{Static}} & \multicolumn{1}{c}{\textbf{Myopic}}
& \multicolumn{1}{c}{\textbf{Index}} & \multicolumn{1}{c}{\textbf{Static}} & \multicolumn{1}{c@{}}{\textbf{Myopic}} \\
\hline
& \multicolumn{3}{c}{(a) $\alpha_k \sim
U[1.05,1.20]$} & \multicolumn{3}{c}{(b) $\alpha_k \sim U[1.20,1.35]$}\\
\textit{MIN} & 0.0000 & 0.0187 & 1.2278 & 0.0000 & 0.0987 & 1.1529 \\
\textit{LQ} & 0.2446 & 4.7749 & 2.4854 & 0.0556 & 8.3715 & 2.6063 \\
\textit{MED} & 0.6471 & 10.9720 & 4.5413 & 0.5215 & 14.7471 & 4.9759 \\
\textit{UQ} & 2.6301 & 17.0301 & 7.0980 & 2.0182 & 21.1644 & 8.8432 \\
\textit{MAX} & 10.8450 & 28.0785 & 22.3554 & 9.5897 & 31.7000 & 22.5440 \\
[4pt]
\textit{N} & \multicolumn{1}{c}{2000} & \multicolumn{1}{c}{2000}
& \multicolumn{1}{c}{2000} & \multicolumn{1}{c}{2000} & \multicolumn{1}{c}{2000}
& \multicolumn{1}{c@{}}{2000} \\
[4pt]
& \multicolumn{3}{c}{(c) $\alpha_k \sim
U[1.35,1.50]$} & \multicolumn{3}{c}{(d) $\alpha_k \sim U[1.50,1.65]$}
\\
[4pt]
\textit{MIN} & 0.0000 & 0.3388 & 1.1130 & 0.0000 & 0.9814 & 0.9718 \\
\textit{LQ} & 0.0122 & 11.2186 & 2.8107 & 0.0034 & 14.3835 & 3.6829 \\
\textit{MED} & 0.2554 & 17.4297 & 5.9093 & 0.1743 & 21.1017 & 7.6089 \\
\textit{UQ} & 1.7601 & 24.0923 & 10.6612 & 1.6311 & 27.6231 & 13.4215 \\
\textit{MAX} & 8.0043 & 33.8457 & 22.7322 & 6.4821 & 36.3746 & 24.4466 \\
[4pt]
\textit{N} & \multicolumn{1}{c}{2000} & \multicolumn{1}{c}{2000} & \multicolumn{1}{c}{2000}
& \multicolumn{1}{c}{2000} & \multicolumn{1}{c}{2000} & \multicolumn{1}{c@{}}{2000}\\
\hline
\end{tabular*}
\end{table}

The results in Table~\ref{tabs5} concern a very simple model in which
the transition rates are state independent. We take $\xi_k(\cdot)
\equiv1,  k = 1,  2$, while the $\eta_k(\cdot)$ also are constant,
with values obtained by sampling from the uniform distribution on
$[0.75,1.25]$. Resource is available to the assets at total rate $R=5$
throughout. In all cases, asset return rates are increasing concave in
the asset state and given by
\[
d(n) = n(n+1)^{-1}, \qquad 0 \leq n \leq10.
\]

These asset management problems prove challenging and the myopic
proposal performs poorly in Table~\ref{tabs5}, being consistently
outperformed by both Index and Static. Over the 2000 problems sampled,
the percentage suboptimality of Index is roughly uniformly distributed
on the interval $[0.0,1.9]$, while that for Static is also roughly
uniform, but across the considerably wider range $[0.0,13.7]$.\vadjust{\goodbreak}

For the next group of problems we set $R=10$ and introduce state
dependence into the transition rates. In~(\ref{s50}) and~(\ref{s51})
we take
%
\begin{equation} \label{s52}
\xi_k(n_k) = \{11^{\alpha_k} - (n_k+1)^{\alpha_k}\}(n_k+1)^{-\alpha
_k+1},\qquad 0 \leq n_k \leq9,
\end{equation}
and
%
\begin{equation} \label{s53}
\eta_k(n_k) = n_k,\qquad 1 \leq n_k \leq10,
\end{equation}
where in~(\ref{s52}) and~(\ref{s53}), $\alpha_k>1$ is a positive
constant. The choices in~(\ref{s52}),~(\ref{s53}) feature in the
numerical study undertaken by Glazebrook, Kirkbride and Ruiz-Hernandez \cite{gla06} of
their much simpler spinning plates model. The function $\xi_k$ in
(\ref{s52}) is decreasing and convex over the range $0 \leq n_k \leq
9$. The
degree of curvature of the function and the value of $\xi_k(0)$ both
increase with the value of $\alpha_k$. For the problems featured in
Table~\ref{tabs6}, we obtain the $\alpha_k$ by sampling from the
uniform distribution on $[1.05,1.50]$. Here the models are such that
achieving improvements to asset performance is increasingly difficult
for higher states. This effect will be most marked when\vadjust{\goodbreak} $\alpha_k$ is
close to the top of its range. Finally, our choice of asset return rate is
%
\begin{equation} \label{s54}
d(n) = \cases{
0, &\quad $0 \leq n \leq4$, \cr
(n-4)/5, &\quad $5 \leq n \leq8$, \cr
1, &\quad $n = 9, 10$.}
\end{equation}
Here state 9 is the minimum for an asset to generate returns at maximal
rate. Further, should an asset deteriorate to the point that its state
is 4 or less it is incapable of generating any returns. In contrast to
the problems featured in Table~\ref{tabs5}, this return is nonconcave
in state.

Please find the results for this group of 2000 problems in Table \ref
{tabs6}. In Table~\ref{tabs7} we consider a slightly modified set of
such problems for which $R=5$ and the downward transition rates are
given by
\[
\eta_k(n_k) = 0.5 n_k,\qquad 1 \leq n_k \leq10.
\]

The problems featured in Tables~\ref{tabs6} and~\ref{tabs7} prove
relatively unchallenging to both Index and Static, in part because of
the highly discrepant upward transition rates obtained from distinct
$\alpha_k$. If we \textit{tame} this feature by rescaling the functions
$\xi_k$ (after $\alpha_k$ has been chosen) such that $\xi_k(0)$ is a
fixed amount (here taken to be 12) then the problems become very much
more difficult and the performance of Static can become quite poor.
Table~\ref{tabs8} features 8000 such problems. The subtables
correspond to distinct ranges for the sampled~$\alpha_k$. In Table
\ref{tabs8}(a)--\ref{tabs8}(d) we have $\alpha_k \sim U[1.05,1.20],
\alpha_k \sim U[1.20,1.35], \alpha_k \sim U[1.35,1.50]$ and $\alpha_k
\sim U[1.50,1.65]$, respectively.\vadjust{\goodbreak} Problem details are otherwise as for
Table~\ref{tabs6}. From Table~\ref{tabs8}, the relatively poor
performance of both Static and Myopic makes it clear that these are
difficult problems for which dynamic policies, which take adequate
account of the future impact of current decisions, really are needed.
The greedy index heuristic delivers a readily understood proposal which
continues to perform robustly even in this very challenging problem
environment. It is especially effective for the problems with larger
sampled $\alpha_k$ in which it is most difficult to maintain assets in
strongly performing states.\vspace*{-3pt}

\section{Conclusions and proposals for further work}\label{sec5}

The paper has described radical extensions to index theory which
facilitate the analysis of dynamic resource allocation problems in
which a single key resource may be assigned more flexibly than is
allowed in classical bandit models. The resulting greedy index
heuristic has been shown to perform strongly for a range of models
which relate to applications, \textit{inter alia}, in queueing control
and asset management which are of independent interest.

Without doubt, the primary obstacle to general implementation of the
approach described concerns the requirement to establish full
indexability. This is that optimal solutions to the single project
problems $P(k,W), 1 \leq k \leq K$, derived from a Lagrangian
relaxation of the original problem, exhibit a property of assigning
diminishing levels of resource uniformly over project states as the
resource charge $W$ increases. While we have been able to demonstrate
this for the models of Sections~\ref{sec3} and~\ref{sec4}, it presents a formidable
challenge in many problems. We propose to develop our approach further
by exploring the quality of index heuristics derived from \textit
{strongly performing} (though possibly not optimal) policies for the
$P(k,W), 1 \leq k \leq K$, which have the above structural property
required to create indices.\vspace*{-3pt}

\section*{Acknowledgment}

We gratefully acknowledge the helpful comments of an anonymous referee
 for challenging us to strengthen the paper.\vspace*{-2pt}


%
\printaddresses

\end{document}